\documentclass[review,onefignum,onetabnum]{siamart171218}

\usepackage{graphicx}
\usepackage{underscore}
\usepackage{cancel}
\usepackage[utf8]{inputenc}
\usepackage{textpos}
\usepackage{lipsum}
\usepackage[]{geometry}     % make the page margins visible
\usepackage[english]{babel}
\usepackage{mathtools}
\usepackage{faktor}
\usepackage{hyperref}
\usepackage{subcaption}
\usepackage{physics}
\usepackage{tikz-cd} 
\usepackage{float}
\usepackage{tikz}
\usepackage{blindtext}
\usepackage{adjustbox}
\usepackage{mathrsfs}
\usepackage{bm}
\usepackage{polski}

\usepackage{multimedia}
\usepackage{hyperref}
\usepackage{wrapfig}
\usepackage{paralist}
\usepackage{bbm}

\usepackage{amsfonts}
\usepackage{epstopdf}
\usepackage{algorithmic}
\usepackage{physics}
\ifpdf
  \DeclareGraphicsExtensions{.eps,.pdf,.png,.jpg}
\else
  \DeclareGraphicsExtensions{.eps}
\fi

\definecolor{KUL_blue}{rgb}{0.11,0.55,0.69}
\definecolor{KUL_dark_blue}{rgb}{0.067, 0.431, 0.541}
\definecolor{KUL_light_blue}{rgb}{0.86, 0.91, 0.94}

\makeatletter
\newcommand\footnoteref[1]{\protected@xdef\@thefnmark{\ref{#1}}\@footnotemark}
\makeatother

\usepackage{tcolorbox}
\newtcolorbox{block}[1]{
    before=\par\smallskip\centering,
	colback=KUL_light_blue,
	colbacktitle=KUL_blue,
	coltitle=white,
	coltext=black,
	colframe=KUL_blue,
	boxrule=1pt,
	arc=2pt,
	fonttitle=\bfseries,
	title={\strut#1}
}

\usepackage{natbib} % Citation styles

\newsiamremark{remark}{Remark}
\newsiamremark{hypothesis}{Hypothesis}
\crefname{hypothesis}{Hypothesis}{Hypotheses}
\newsiamthm{claim}{Claim}

\headers{Beyond Classical Derivatives}{C. Verbeeck, N. Sfakianakis}

\title{Beyond Classical Diffusion: Fractional Derivatives in Transport and Stochastic Systems\thanks{Submitted to the editors \today.
}}

\author{Cypres Verbeeck\thanks{Mathematical Biology Group, University of St Andrews, Scotland.}
\and Nikolaos Sfakianakis\footnotemark[2]}

\usepackage{amsopn}

\makeatletter
\newcommand*{\addFileDependency}[1]{% argument=file name and extension
  \typeout{(#1)}% latexmk will find this if $recorder=0 (however, in that case, it will ignore #1 if it is a .aux or .pdf file etc and it exists! if it doesn't exist, it will appear in the list of dependents regardless)
  \@addtofilelist{#1}% if you want it to appear in \listfiles, not really necessary and latexmk doesn't use this
  \IfFileExists{#1}{}{\typeout{No file #1.}}% latexmk will find this message if #1 doesn't exist (yet)
}
\makeatother

\newcommand*{\myexternaldocument}[1]{%
    \externaldocument{#1}%
    \addFileDependency{#1.tex}%
    \addFileDependency{#1.aux}%
}
\ifpdf
\hypersetup{
  pdftitle={Beyond Classical Derivatives},
  pdfauthor={C. Verbeeck, N. Sfakianakis}
}
\fi

\myexternaldocument{ex_supplement}

\usepackage[normalem]{ulem}

\newcommand{\nkschange}[2]{%
  \ifmmode
    {\color{red}\cancel{#1}}{\color{blue} ns:#2}
  \else
    {\color{red}\sout{#1}}{\color{blue} ns:#2}
  \fi
}
\newcommand{\R}{\mathbb{R}}

\usepackage{tcolorbox}
\usepackage{hyperref}
\newtcolorbox{blocktitle}[1][]{ % 1 optional argument for the title
  colback=white, % Background color
  colframe=white, % Frame color (invisible)
  sharp corners, % No rounded corners
  boxrule=0pt, % No border
  fonttitle=\bfseries, % Bold title
  title=#1, % Title text
  left=10pt, % Left indent
}

\begin{document}

% \nkscomm{Possible journals: Multiscale modelling an simulations}

% \nkscomm{Possible editor: Jose-Antonio Carrillo}

% \nkscomm{Possible reviewers: 1) Jose's co-author in the fractional drv paper, 2) Franz Achleitner TU  Wien, 3) Dietmar Oelz Queensland university}

\maketitle

% REQUIRED
\begin{abstract}
% Differential operators of integer order are known to capture local and isotropic effects, both in space and time. However, especially in biology, the increasing need to model complex phenomena with underlying properties such as spatial heterogeneity requires new modelling tools. The fractional calculus framework enables the development of more sophisticated models that capture the complex dynamics inherent to various biological systems. This talk will focus on how fractional reaction-diffusion equations naturally arise and can be derived from continuous-time random walks, highlighting the role of heavy-tailed distributions in the process. Both fractional partial differential equations, on the macroscopic level, as well as fractional stochastic differential equations, on the microscopic level, will be examined. For simple Riesz-fractional diffusion models, we will showcase comparative simulations, highlighting the key differences between fractional and classical diffusion.
% We propose a new numerical scheme that implements periodic boundary conditions, to control the loss of mass density. 
Integer-order differential operators were originally used to describe local and isotropic effects, in both space and time. However, in fields like biology, the modelling of complex phenomena with spatial heterogeneity necessitates more advanced approaches. The fractional calculus framework provides powerful tools for developing models that better capture the intricate dynamics of biological systems. This paper derives fractional reaction-diffusion equations from continuous-time random walks, highlighting the role of heavy-tailed distributions in the process. Both fractional partial differential equations, on the macroscopic level, as well as fractional stochastic differential equations, on the microscopic level, will be derived and simulated from, for simple Riesz-fractional diffusion models. A new numerical scheme that implements periodic boundary conditions is proposed to control the loss of mass density. We highlight the key differences between fractional and classical diffusion.
\end{abstract}

% REQUIRED

% REQUIRED

\section{Introduction}
Differential operators of integer order are known to capture local and isotropic effects, both in space and time. However, especially in biology, the increasing need to model complex phenomena with underlying properties such as spatial heterogeneity and/or effects of memory requires new (and more complex) modelling tools, \cite{Motivation, West2016}. The fractional calculus framework enables the development of more sophisticated models that accurately capture the complex dynamics inherent in various systems: including viscoelasticity, \cite{Mainardi2010}, hydrological processes, \cite{Benson2000}, and biological applications such as bioengineering, \cite{Magin2006}. Fractional calculus could even be applied to describe motor control systems, notably the dynamics of motor and premotor neurons, \cite{Anastasio1994}.
%This article extends the use of fractional calculus to describe the dynamics of motor and premotor neurons. It suggests that the oculomotor integrator may be of fractional order. The fractional derivative dynamics of motor and premotor neurons may serve to compensate fractional integral dynamics of the eye. Fractional differentiation can be used to account for the constant phase shift across frequencies, and the apparent decrease in time constant as VOR and pursuit frequency increases, that are observed for motor and premotor neurons. Fractional integration can reproduce the time course of motor and premotor neuron saccade-related activity, and the complex dynamics of the eye. 
%Fractional dynamics may be applicable not only to the oculomotor system, but to motor control systems in general.
%The fractional-order dynamics of brainstem vestibulo-oculomotor neurons

Specifically, we want to see for ourselves how fractional reaction-diffusion equations arise in a natural way, starting from basic principles. Similar to how one would derive the macroscopic diffusion equation, a continuous-time random walk (CTRW), together with mass conservation equations, provide the starting point of our derivations. The use of CTRWs in fractional diffusion modelling has been extensively studied, particularly in the context of numerical simulations of individual particle trajectories, by \cite{Kleinhans2007}. \cite{SubdiffFD} approached fractional diffusion in 2004, from a probabilistic perspective. Inspired by their methodology and results, this paper transfers their ideas to the slightly different density setting, where creation-destruction phenomena are also factored into the system's progression. Starting from the mesoscopic density equations provided by \cite{reaction-transport}, we are able to argue how fractional derivatives naturally make their appearance in the resulting macroscopic description of the system. 

Diffusion describes the physical process of multiple particles undergoing a random walk on a microscopic level. Integer-order diffusion equations present the most well-known elementary models to describe the diffusion process. They assume that the stochastic walks follow a standard Gaussian distribution. However, this simplistic underlying assumption neglects to incorporate lesser effects of complicated biological scenarios, which are nonetheless equally important as more prominent (collective) movements. Moreover, there are many diffusion processes in nature that do not satisfy the Gaussianity assumption. Hence, employing the integer-order diffusion equation does not always lead to an accurate representation of diffusion as observed in reality. In porous media, for instance, anomalous diffusion is perceived to be dominant. To quote \cite{XieFang}: ``A large number of experiments show that the fractional diffusion equations can provide a more accurate mathematical description of these anomalous diffusion processes in nature than the integer-order diffusion model [...]''. 
%\nkscomm{A direct quote works better only if the author is extremely well-known, like Pythagoras, Galois, Einstein, Nikos}
For a detailed review of anomalous diffusion and fractional dynamics, see also \cite{Metzler2000} or \cite{Henry2000}. For a broader perspective on fractional dissipative PDEs---including reaction-diffusion systems, their stability, bifurcation and travelling waves analysis, as well as some numerical simulations---we refer to \cite{Achleitner, FracTheory}.

Fractional diffusion equations are starting to attract more and more interest. Notably, Riesz-fractional diffusion equations have been adopted in many scientific and engineering applications. Often, an exact analytical solution of one such model is either not available or too complicated to reproduce. This is why numerical methods for the simulation of fractional (diffusion) models need to be developed. We refer the interested reader to Li and Zeng's \textit{Numerical Methods for Fractional Calculus} for an overview of numerical methods for the computation of various fractional differential equations, \cite{Li2015}. Recent advances have introduced memory-efficient spectral methods for solving time-fractional PDEs, significantly reducing computational costs while maintaining accuracy, see \cite{Carrillo2023}. Although various numerical methods have already been introduced for Riesz (one-dimensional) fractional diffusion in particular, \cite{XieFang}, many schemes fail to correctly account for mass conversation. 
%The fractional Laplacian (−Δ)α/2 
%A numerical method for the fractional Laplacian is proposed, based on the singular integral representation for the operator. The method combines finite difference with numerical quadrature, to obtain a discrete convolution operator with positive weights.  %treatment of far field boundary conditions using an asymptotic approximation to the integral is used to obtain an accurate method. 
Most available numerical methods assume Dirichlet boundary conditions.
%study of solution techniques for problems involving fractional powers of symmetric coercive elliptic operators in a bounded domain with Dirichlet boundary conditions. we propose a truncation that is suitable for numerical approximation. We discretize this truncation using first degree tensor product finite elements. We derive a priori error estimates in weighted Sobolev spaces. 
See also \cite{Nochetto2015} for a study of solution techniques for problems involving fractional powers of symmetric coercive elliptic operators in a bounded domain with Dirichlet boundary conditions. In \cite{Huang2014}, far-field boundary conditions are assumed to solve the fractional Laplacian numerically. Hermite spectral methods lie at the heart of solving fractional PDEs in unbounded domains in \cite{Mao2017}, while a spectral-Galerkin method solves ``integral fractional Laplacian PDEs'' in unbounded domains of $\mathbb{R}^d$ in \cite{Sheng2020}.
Nevertheless, plenty of real-life problems involving space-fractional diffusion on a bounded domain require mass-conserving boundary conditions instead, \cite{BCtwosided}. 
%\cite{Cusimano2018}

This paper aims to compare fractional against regular diffusion on both the microscopic as well as the macroscopic level. The microscopic model calls for stochastic differential equations, drawing from normal and stable distributions. On the macroscopic level, we propose a new numerical scheme that incorporates periodic boundary conditions, by closing up a one-dimensional spatial domain. By doing so, we avoid any blow-up phenomena that may occur at the boundary of the domain, see also \cite{Abatangelo2023} for an elaborate discussion of the explosive behaviour characteristic of equations driven by non-local---fractional---operators.
%Here, we show that, for equations driven by a wide class of nonlocal fractional operators, different blow-up phenomena may occur at the boundary of the domain. We describe such explosive behaviours
Existing approximations will form the basis of a modified Gr{\"u}nwald-Letnikov discretisation for the spatial Riesz derivative. Another possible discretisation of the spectral fractional Laplacian on a bounded domain is discussed by, for instance, \cite{Cusimano2018}. Their formalism is naturally made to handle different types of boundary constraints.
%In this work, we propose novel discretizations of the spectral fractional Laplacian on bounded domains based on the integral formulation of the operator via the heat-semigroup formalism.  can be implemented on possibly irregular bounded domains, and can naturally handle different types of boundary constraints. 
We use the Gr{\"u}nwald-Letnikov operator approximation at the midpoint of the domain and replicate the resulting stencil at every other lattice point. We will do so by introducing a fixed vector carrying the weights of the sum, all the while permuting the position indices of all spatial grid points, to account for ghost points and enforce symmetry of the scheme itself. The periodicity acts to transfer information across boundaries and was introduced to control the mass in finite domains. 
%\nkschange{Direct comparison of the dynamics resulting from an---inherently---fractional to a regular model indicates greater \nkscomm{mass} loss for fractional diffusion, by a factor of order $10^{7}$.}{} 

Whilst succeeding in controlling the loss of mass at the endpoints, a complete conservation of mass is not observed, for neither regular nor fractional (of order $1.5$) diffusion. Of course, the choice of different parameters, such as final time $T$, temporal step size $\tau$, and size of the domain itself, affect the results quantitatively. Larger $T$ and narrower domains $(a,b)$---independently---are responsible for more significant loss rates. A smaller time step $\tau$ or a bigger domain acts to decrease the total loss of mass. After all, in the latter case, effects at the boundary become more and more negligible and boundary conditions can be avoided altogether. 

The remainder of this paper is organised as follows. In Section \ref{sec:derivation}, we derive the fractional reaction-diffusion equations from CTRWs, highlighting the role of heavy-tailed distributions in the derivation of fractional derivatives. Section \ref{sec:fractSDE} explores the connection between microscopic stochastic processes and macroscopic fractional diffusion, introducing fractional stochastic differential equations (FSDEs). In Section \ref{sec:numerics}, we provide comparative simulations that exhibit the difference between fractional and classical diffusion, both at the microscopic and macroscopic level.

%close the lattice (ghost values)
%BCs that keep \nkscomm{mass} inside (so it doesn't flow out of the domain)
%enforce conservation in finite domain somehow 
%introduce ghost cells for \nkscomm{mass} conservation to either side of both endpoints

%================================================================================
\section{Derivation of a Fractional Reaction-Diffusion Density Equation}\label{sec:derivation}
The connection between CTRWs and fractional diffusion has been extensively studied before, in \cite{Metzler2000}, where they put the focus on subdiffusive processes. 
Starting from the setup of a CTRW, we will first derive diffusion equations for the probability $p(x,t)$ of an agent being at a position $x$ at time $t$. The agent jumps around in a one-dimensional domain. The law of total probability will enable us to write down the governing equations for the random walker. In our quest for a more explicit expression for $p$ in terms of (presumably) known quantities, namely the waiting-time and jump-length distributions, we have to resort to Fourier-Laplace integral transforms. The first intermediate result is obtained by working out the Fourier-Laplace transform of $p$, which leads to the \textit{Montroll-Weiss equation}.

Since we are, moreover, interested in the macroscopic behaviour of the particle, we have to take suitable limits (corresponding to the correct space-time scaling) in the latter equation. Different asymptotics of the waiting time- and probability distributions (notably, the heavy-tailedness thereof) result in quantitatively different expressions for the Fourier-Laplace transform of $p$. This probability transform can then be further simplified by adopting the governing scaling relation. Upon translating the obtained results in terms of the probability $p$ again, by taking the inverse Fourier-Laplace transforms, we arrive at a (fractional) diffusion equation that is weakly asymptotically fulfilled by $p$. There are two types of fractional derivative formalisms that arise naturally in this setting.

Repeating the same process, but now adopting a density framework, we can similarly arrive at general fractional reaction-diffusion equations satisfied by the particle density. This time, one should start with mean-field equations for the mesoscopic density. A similar process, guided by Fourier-Laplace transformations (and their inverses)---under suitable limits---lets us incorporate fractional orders in the derivatives modelling reaction-diffusion.

\subsection{Continuous-Time Random Walks}\label{CTRWs} 
Suppose that we let a single particle undergo a one-dimensional CTRW. This means that after a specific (waiting) time $t_i$, the particle will have moved over a certain distance $x_i$ (note that this can be a non-positive quantity, in general). Remark that both the waiting times sequence $(\tau_i)_i\geq 0$ and the jumps of the diffusing particle $(x_i)_i$ are independent and identically distributed, according to the probability density functions (PDFs) $\eta(t)$, $w(x)$. In general, the density $\phi(x,t)$ represents the joint waiting-time- and jump length distribution. We can, as such, view the waiting times as representing a point process, with the jumps being the associated ``rewards'', \cite{reaction-transport}. %describes jump process subordinate to a renewal process

Writing
\begin{equation}
    t_n = t_0 + \sum_{i=1}^n \tau_i, \hspace*{0.2cm} (t_k - t_{k-1} = \tau_k \text{ for } k\geq 1)
\end{equation}
for the total time passed upon occurrence of the $n$-th jump exactly. Here, the $\tau_i$ represent the individual waiting times between consecutive shifts in the position. 
The position $x(t)$ of the particle at time $t$ is then recovered by adding the individual distances $x_i$, travelled up until that time $t$, to the initial position $x(0)=x_0$. 
\begin{equation}
    x(t) = x_0 + \sum_{i=1}^{N(t)} x_i.
\end{equation}
Note that here, we also need to take into account the index $N(t)$ which is defined as the exact number of jumps taking place during the time interval $[0,t]$, namely (\cite{reaction-transport})
\begin{equation}
    N(t) = \max_{n} \{ t_n \leq t\}.\label{Nt}
\end{equation}

\subsection{Derivation of the Montroll-Weiss Equation}
%\nkscomm{As this derivation is known, and since we want to save space, I would suggest that we move this sub-section to the Appendix and include here only the basic results/outcomes. We can then refer to the Appendix and the Mendez eta al}
In this section, we derive the Montroll-Weiss equation for a single particlecontinuous-time random walk. The method of proof is based on the derivation described in \cite{reaction-transport}, p.61 onward. The specifics are discussed here in greater detail.

The setup can be described as follows, let us consider a single particle undergoing a CTRW in one dimension. Using a probabilistic approach, we can describe the particle's position at every time $t$. 
Starting from the law of total probability, one may write
\begin{align}
    \underbrace{p(x,t)}_{\substack{\text{probability that a particle is}\\\ \text{at position $x$ at time $t$}}} &= \underbrace{p_0(x)}_{\substack{\text{probability that a particle is}\\ \text{at position $x$ at time $t=0$}\\ \text{(initial PDF)}}}\overbrace{\Psi(t)}^{\text{``survival probability''}} 
    \nonumber \\
    &
    +\int_0^t\underbrace{j(x,t-\tau)}_{\substack{\text{probability that a particle}\\ \text{reaches point $x$ at}\\ \text{time $t-\tau$ exactly}}}\Psi(\tau)d\tau, \label{massprob}
\end{align}
where 
\begin{align}
    \Psi(t) &= \int_t^{\infty}\int_{\R} \phi(y,\tau)dyd\tau 
    = 1-\int_0^{t}\int_{\R} \phi(y,\tau)dyd\tau \label{waitingtime}  
\end{align}
denotes the probability that the particle doesn't jump at all after a waiting time less than $t$ (i.e. the probability of no change in the particle's position), \cite{reaction-transport}. 
Moreover, $j(x,t)$ can further be decomposed as follows 
\begin{align}
    \underbrace{j(x,t)}_{\substack{\text{probability that the particle reaches}\\ \text{position $x$ by time $t$}\\ \text{(by means of transport)}}} &= \int_{\R} p_0(x-z)\underbrace{\phi(z,t)}_{\substack{\text{probability that the particle has moved}\\ \text{a distance $z$ after a waiting time $t$}}}dz \nonumber\\
    & + \underbrace{\int_0^t\int_{\R} j(x-z,t-\tau)\phi(z,\tau)dzd\tau}_{\substack{\text{probability that the particle is situated at the position $x-z$}\\ \text{at time $t-\tau$, and has moved over a distance $z$ after a waiting time $\tau$}}}. \label{jprob}
\end{align}

We wish to solve equation \eqref{massprob} for $p(x,t)$ explicitly. However, due to its integral dependence on the probability $j(x,t)$, which is implicitly defined in terms of itself through \eqref{jprob}, we will have to resort to certain integral transforms of the equations. Specifically, we will take Fourier-Laplace transforms of the above equations and ultimately arrive at the Montroll-Weiss equation. 

First of all, taking the Fourier transform $\mathcal{F}$---with respect to the space variable $x$---of \eqref{massprob} leads us to
\begin{align}
    \mathcal{F}[p](k,t) &= \mathcal{F}[p_0](k)\Psi(t)+\int_0^t\mathcal{F}[j](k,t-\tau)\Psi(\tau)d\tau=\mathcal{F}[p_0](k)\Psi(t)+\mathcal{F}[j](k,\cdot)*_t\Psi(\cdot), \label{fourierp}
\end{align}
where `$*_t$' was used to denote the convolution between functions with respect to the time variable.
Here, it should be noted that the operator $\mathcal{F}$ acting on a function $f$ is assumed to return the angular frequency, non-unitary Fourier transform of $f$. The variable $k$ is used to denote the angular frequency.
Consecutively taking the Laplace transform of \eqref{fourierp} (moving from the time variable $t$ to the new variable $s$ of the transform, keeping in mind the multiplicative property of the aforementioned transform with respect to convolution products) then brings us to 
\begin{align}
    \mathcal{L}[\mathcal{F}[p](k,\cdot)](s)
    &=\mathcal{F}[p_0](k)\mathcal{L}[\Psi](s)+\mathcal{L}[\mathcal{F}[j](k,\cdot)*_t\Psi(\cdot)](s) \nonumber\\
    &=\mathcal{F}[p_0](k)\mathcal{L}[\Psi](s)+\mathcal{L}[\mathcal{F}[j](k,\cdot)](s)\mathcal{L}[\Psi(\cdot)](s).\label{FL1}
\end{align}
We now aim to find explicit expressions for both $\mathcal{L}[\Psi](s)$ and $\mathcal{L}[\mathcal{F}[j](k,\cdot)](s)$. Upon taking the Fourier transform of \eqref{jprob}, by definition
\begin{align}
    \mathcal{F}[j(\cdot, t)](k)&=\int_{\R}j(x,t)e^{-ikx}dx\nonumber\\
    &= \underbrace{\int_{\R}\int_{\R}p_0(x-z)\phi(z,t)e^{-ikx}dzdx}_{I_1}
    + \underbrace{\int_{\R}e^{-ikx}dx \int_0^t\int_{\R}j(x-z, t-\tau)\phi(z,\tau)dzd\tau}_{I_2}.\label{I1plusI2}
\end{align}
Assuming that we may apply Fubini 
%{\color{red}make this more precise?} 
to change the order of integration in the integrals above, one finds 
\begin{align}
   I_1 &=\int_{\R}\phi(z,t)\int_{\R}p_0(x-z)e^{-ik(x-z)}e^{-ikz}dxdz=\int_{\R}\phi(z,t)e^{-ikz}\int_{\R}p_0(u)e^{-iku}dudz \label{changeofvar}\\
   &=\int_{\R}\phi(z,t)e^{-ikz}\mathcal{F}[p_0](k)dz=\mathcal{F}[\phi(\cdot,t)](k)\mathcal{F}[p_0](k),\label{fourierI1}
\end{align}
where we have performed the change of variables $u = x-z$ (for $z=z$ unmodified) in \eqref{changeofvar} which has corresponding Jacobian
$\begin{pmatrix}
    1 & -1 \\ 0 & 1
\end{pmatrix}.$
Similarly, the second integral reduces to 
\begin{align}
    I_2&=\int_0^t d\tau \int_{\R} \phi(z,\tau)dz \int_{\R}e^{-ikx}j(x-z,t-\tau)dx=\int_0^t d\tau \int_{\R} \phi(z,\tau) e^{-ikz}\mathcal{F}[j(\cdot,t-\tau)](k)dz\nonumber\\
    &=\int_0^t \mathcal{F}[\phi(\cdot,\tau)](k)\mathcal{F}[j(\cdot,t-\tau)](k) d\tau= \mathcal{F}[\phi](k,\cdot)*_t\mathcal{F}[j](k,\cdot).\label{fourierI2}
\end{align}
Together, \eqref{fourierI1} and \eqref{fourierI2} imply, by \eqref{I1plusI2}, that
\begin{align}
    \mathcal{F}[j](k,t) &= \mathcal{F}[\phi](k,t)\mathcal{F}[p_0](k) + \mathcal{F}[j](k,\cdot)*_t\mathcal{F}[\phi](k,\cdot). \label{fourierjtotal}
\end{align}
Hence, consecutively taking the Laplace transform of the expression for $\mathcal{F}[j](k,t)$ above, one arrives at
\begin{align}
    \mathcal{L}[\mathcal{F}[j]](k,s) &= \mathcal{L}[\mathcal{F}[\phi]](k,s)\mathcal{F}[p_0](k) + \mathcal{L}[\mathcal{F}[j]](k,s)\mathcal{L}[\mathcal{F}[\phi]](k,s)
\end{align}
due to linearity. 
Therefore, upon rearranging, the equality
\begin{align}
    \mathcal{L}[\mathcal{F}[j]](k,s) &=\frac{\mathcal{L}[\mathcal{F}[\phi]](k,s)\mathcal{F}[p_0](k)}{1-\mathcal{L}[\mathcal{F}[\phi]](k,s)} \label{laplacefourierj}
\end{align}
rolls out altogether. As such, we have found an expression for $\mathcal{L}[\mathcal{F}[j]](k,s)$ in terms of the---presumed priorly known---densities $p_0$ and $\phi$. 

As far as $\Psi(t)$ goes, we continue with the Laplace transform of \eqref{waitingtime}. By definition of the Laplace transform, one finds
\begin{align}
    \mathcal{L}[\Psi](s) &= \int_0^{\infty} \Psi(t) e^{-st}dt
    = \int_0^{\infty} \int_t^{\infty}\int_{\R} \phi(y,\tau)e^{-st}dyd\tau dt= \int_{\R} dy \int_0^{\infty}d\tau \int_0^{\tau} \phi(y,\tau)e^{-st}dt. \label{laplacepsiint}
\end{align}
Note that we have again appealed to Fubini's Theorem to modify the order of integration (accounting for the bounds) in this triple integral.
%(where one has to act carefully as to change the bounds alongside the order, where $0\leq t\leq \tau <\infty$ is assumed).  

Explicit computation of the integrals in \eqref{laplacepsiint} then leads to the following ($s\neq0$)
\begin{align}
    \mathcal{L}[\Psi](s)
    &= \int_{\R} dy \int_0^{\infty}d\tau  \phi(y,\tau)\left[-\frac{e^{-st}}{s}\right]_{t=0}^{t=\tau}= \int_{\R} dy \int_0^{\infty}d\tau \phi(y,\tau)\left[\frac{1-e^{-s\tau}}{s}\right]\nonumber\\
    &= \frac{1}{s}\left\{\int_{\R} dy \int_0^{\infty}d\tau \phi(y,\tau)-\int_{\R} dy \int_0^{\infty}d\tau \phi(y,\tau)e^{-s\tau}\right\}\nonumber\\
    &= \frac{1}{s}\left\{1-\int_{\R} dy\mathcal{L}[\phi(y,\cdot)](s)\right\},\label{laplacepsiint2}
\end{align}
again upon straightforward application of the definition of the Laplace transform of a function. Note also that $\phi(y,t)$ represents a probability density and therefore integrates to one over the space- and time domain.

\textbf{Independence Assumption.} For simplicity and ease of calculation, we have assumed the jump lengths and waiting times to be independent of one another. As a result of this assumption, we factorise the joint PDF $\phi(y,t)=w(y)\eta(t)$ in terms of the individual waiting time probability density function $\eta(t)$ and the jump length probability density function $w(y)$ (also known as the dispersal kernel, in the uncoupled case). 

Using the independence assumption in our calculations \eqref{laplacepsiint2}, we can further simplify its last integral. With $\mathcal{L}[\phi(y,\cdot)](s)$ now equal to $\mathcal{L}[\eta](s)w(y)$, altogether
\begin{align}
    \mathcal{L}[\Psi](s) &= \frac{1}{s}\left\{1-\mathcal{L}[\eta](s)\int_{\R} dy  w(y)\right\}%\nonumber\\
    = \frac{1}{s}\left\{1-\mathcal{L}[\eta](s)\right\}\label{laplacepsifinal}
\end{align}
(with $w(y)$ is a probability density in its own regard). 
Combining \eqref{laplacefourierj} and \eqref{laplacepsifinal} in \eqref{FL1}, we arrive at the so-called \textit{Montroll-Weiss equation}. 
\begin{align}
    \mathcal{L}[\mathcal{F}[p](k,\cdot)](s)
    &=\left\{\mathcal{F}[p_0](k)+\mathcal{L}[\mathcal{F}[j](k,\cdot)](s)\right\}\mathcal{L}[\Psi(\cdot)](s)\nonumber\\
    &=\left\{\mathcal{F}[p_0](k)+\frac{\mathcal{L}[\mathcal{F}[\phi]](k,s)\mathcal{F}[p_0](k)}{1-\mathcal{L}[\mathcal{F}[\phi]](k,s)}\right\}\frac{1-\mathcal{L}[\eta](s)}{s}\nonumber\\
    &=\frac{\mathcal{F}[p_0](k)}{s}\left\{\frac{1-\mathcal{L}[\mathcal{F}[\phi]](k,s)+\mathcal{L}[\mathcal{F}[\phi]](k,s)}{1-\mathcal{L}[\mathcal{F}[\phi]](k,s)}\right\}\left\{1-\mathcal{L}[\eta](s)\right\}\nonumber\\
    &=\frac{\mathcal{F}[p_0](k)}{s}\left\{\frac{1-\mathcal{L}[\eta](s)}{1-\mathcal{L}[\mathcal{F}[\phi]](k,s)}\right\} \label{MontrollWeiss}
\end{align}

\subsection{Regular Diffusion as a Limiting Case}
We are interested in the macroscopic behaviour of an agent undergoing a CTRW. Hereto, let us consider the large-scale, long-time limit in the last section. We will come to see that this---asymptotically---gives rise to a diffusion equation for the probability $p(x,t)$. Depending on the initial assumptions made regarding waiting time and jump length distributions, our diffusion model will either take on the form of a regular diffusion PDE, or it will be some fractional (anomalous diffusion) counterpart thereof. In this section, we assume that the PDFs are not heavy-tailed, this will result in the familiar Brownian diffusion limiting case.

Remark that we still consider the jump length and waiting time PDFs to be independent, notably $\phi(x,t)=w(x)\eta(t)$. Using the Montroll-Weiss equation \eqref{MontrollWeiss} as our starting point, we first need to translate what ``long-time, large-scale'' means for the variables $s$ and $k$ of the PDFs' integral transforms. As it turns out, the corresponding limits for $k$ and $s$ are given by $k,s \to 0$, see also Section \ref{subdifflim} for the general justification of this fact.

At this point, we also have to enforce some additional assumptions. 
In general, CTRWs can be characterised by their mean waiting time 
\begin{equation*}
    \overline{T} = \int_0^{\infty} t \eta(t)dt,
\end{equation*}
along with the second moment of their jump length distribution
\begin{equation*}
    \sigma^2 = \int_{\R}z^2w(z)dz. 
\end{equation*}
If the waiting time PDF is not heavy-tailed, then the mean waiting time $\overline{T}$ is finite. Moreover, assuming the CTRW to have a short-range jump length PDF $w(x)$, the variance $\sigma^2$ will be finite as well. Assuming an even distribution for the jump lengths, the mean jump length is reduced to $0$. Under all of these assumptions, we are able to Taylor expand both occurring exponentials in the Laplace and Fourier transform definitions. More precisely,
\begin{align*}
    \lim_{s\to 0} \mathcal{L}[\eta](s) &= \lim_{s\to 0} \int_0^{\infty} \eta(t)e^{-st}dt= \lim_{s\to 0} \int_0^{\infty} \eta(t)\left[1-st +\frac{(st)^2}{2} + \mathcal{O}(s^3)\right]dt \hspace*{0.2cm} \\
    &= 1-s\overline{T} + \mathcal{O}(s^2) \hspace*{0.2cm} \text{as $s\to 0$},
\end{align*}
where $\eta$ is assumed to be a probability density. Analogously, one computes
\begin{align*}
    \lim_{k\to 0} \mathcal{F}[w](k) &= \lim_{k\to 0} \int_{\R} w(x)e^{-ikx}dx= \lim_{k\to 0} \int_{\R} w(x)\left[1-ikx-\frac{(kx)^2}{2} + \mathcal{O}(k^3)\right]dx \hspace*{0.2cm} \\
    &= 1-\frac{k^2}{2}\sigma^2 + \mathcal{O}(k^3) \hspace*{0.2cm} \text{as $k\to 0$}.
\end{align*}

Coming back to \eqref{MontrollWeiss}, the Fourier-Laplace (F-L) transform of the probability density $p$ will be given by 
\begin{align}
    \mathcal{L}[\mathcal{F}[p](k,\cdot)](s)&=\frac{\mathcal{F}[p_0](k)}{s}\left\{\frac{1-\mathcal{L}[\eta](s)}{1-\mathcal{L}[\mathcal{F}[\phi]](k,s)}\right\}
    =\frac{\mathcal{F}[p_0](k)}{s}\left\{\frac{1-\mathcal{L}[\eta](s)}{1-\mathcal{F}[w](k)\mathcal{L}[\eta](s)}\right\}\label{LFp}\\
    % &=\frac{\mathcal{F}[p_0](k)}{s}\left\{\frac{1-\mathcal{L}[\eta](s)}{1-\mathcal{F}[w](k)\mathcal{L}[\eta](s)}\right\}\label{LFp}\\
    &=\frac{\mathcal{F}[p_0](k)}{s}\left\{\frac{1-\{1-s\overline{T} + \mathcal{O}(s^2)\}}{1-\{1-\frac{k^2}{2}\sigma^2 + \mathcal{O}(k^3)\}\{1-s\overline{T} + \mathcal{O}(s^2)\}}\right\} \hspace*{0.2cm} \text{as $k,s \to 0$}\nonumber\\
 &\displaybreak[1] 
 \nonumber\\
    &\approx\frac{\mathcal{F}[p_0](k)}{s}\left\{\frac{s\overline{T}}{\frac{k^2}{2}\sigma^2 +s\overline{T}}\right\} \hspace*{0.2cm} \text{as $k,s \to 0$}\nonumber\\
    &\approx\mathcal{F}[p_0](k)\left\{\frac{1}{\frac{k^2\sigma^2}{\overline{2T}} +s}\right\} \hspace*{0.2cm} \text{as $k,s \to 0$}\nonumber\\
    &=\mathcal{F}[p_0](k)\left\{\frac{1}{Dk^2 +s}\right\} \hspace*{0.2cm} \text{as $k,s \to 0$}, \label{LPdif}
\end{align}
where $D=\frac{\sigma^2}{2\overline{T}}$ equals the diffusion coefficient. Rewritten, the Montroll-Weiss equation \eqref{LFp} in the limit scenario $k,s \to 0$ \eqref{LPdif} becomes equivalent to the (approximate) equality
\begin{align}
     (Dk^2 +s)\mathcal{L}[\mathcal{F}[p](k,\cdot)](s)&\approx\mathcal{F}[p_0](k). 
\end{align}

We can now make use of two well-known identities of the Fourier and Laplace integral transforms, namely
\begin{align*}
    \mathcal{F}\left[\pdv[2]{}{x}p(x,t)\right] = -k^2 \mathcal{F}[p(\cdot,t)](k) ~\ \text{ while } ~\
    \mathcal{L}\left[\pdv[]{}{t}p(x,t)\right] = s \mathcal{L}[p(x,\cdot)](s) - p_0(x).
\end{align*}
Inserting these identities in \eqref{LPdif}, we can rewrite the latter equation in the form 
\begin{align}
    \mathcal{F}[p_0](k)&\approx(Dk^2 +s)\mathcal{L}[\mathcal{F}[p](k,\cdot)](s)
    \\ \iff \mathcal{F}[p_0](k) &\approx -D \mathcal{L}\left[\mathcal{F}\left[\pdv[2]{}{x}p(x,t)\right]\right](k,s) + \mathcal{L}\left[\pdv[]{}{t}\mathcal{F}[p(x,t)]\right](k,s) + \mathcal{F}[p_0](k)\label{51}
\end{align}
as $k,s\rightarrow 0$.
We then make use of the fact that the Fourier transform is taken with regard to the spatial variable $x$ (thereby treating $t$ as a fixed constant). As such, we are able to pass the partial derivative with respect to time inside the Fourier operator in \eqref{51}, which results in
\begin{align}
\mathcal{L}\left[\mathcal{F}\left[D\pdv[2]{}{x}p(x,t)\right]\right](k,s) &\approx \mathcal{L}\left[\mathcal{F}\left[\pdv[]{}{t}p(x,t)\right]\right](k,s) \label{52}
    \\ \iff D\pdv[2]{}{x}p(x,t) &\approx \pdv[]{}{t}p(x,t)\label{regdiff}
\end{align}
as $k,s\rightarrow 0$.
Upon taking the Laplace- and then the Fourier inverse in \eqref{52}, we find that the probability density $p(x,t)$ obeys the regular diffusion equation \eqref{regdiff}---in the large-scale, long-time limit that is, assuming that both the waiting-time PDF as well as the jump length PDF do not possess heavy tails, \cite{reaction-transport}.   

\subsection{Subdiffusive Limit}\label{subdifflim}
In the last paragraph, we were able to derive the regular diffusion equation as a limiting case of Montroll-Weiss, under specific conditions. Crucial to the former derivations were the assumptions that the jump length and waiting time PDFs are not heavy-tailed---having finite mean waiting time $\overline{T}$, and variance $\sigma^2$ respectively. In reality, it may be the case that heavy-tailed distributions are observed as far as the agent's spatial-temporal distribution goes. In this case, the Mittag-Leffler function $\Psi(t)=E_{\beta}((-t)^{\beta})$ (where $0<\beta<1$) might be a better fit for the survival probability. It generalises the exponential function, which coincides with the Mittag-Leffler function for $\beta =1$. 

Let us now adopt the assumption that the agent undergoing a one-dimensional CTRW is described by heavy-tailed spatial-and temporal distributions. 
Asymptotically, suppose $\eta(t)\sim t^{-(1+\beta)}$ as $t\rightarrow \infty$ while $w(x)\sim x^{-(1+\alpha)}$ as $x\rightarrow \infty$, where $0<\beta<1$ and with $w(x)\sim x^{-(1+\alpha)}$, \cite{SubdiffFD}. Investigating the long-time, large-scale limit, we recover a \textit{fractional diffusion equation}. The latter describes the sub-(or super-) diffusive behaviour characteristic to the particle. The main theoretical aim of this paper is provided within this section, where it is shown that the solution $p(x,t)$ of the CTRW master equation resulting from \eqref{MontrollWeiss} weakly converges to the solution of a Cauchy problem for the fractional diffusion equation.
%(which has power-law waiting time and jump length distributions)

Proper scaling of the CTRW will lead to a fractional diffusion equation after taking the diffusive limit. This resulting equation generalises the idea of \textit{anomalous} diffusion for power-law tailed probability densities. 
We will still adopt the underlying assumption that jumps and waiting times are not correlated. That is to say, we restrict ourselves to the uncoupled case where jump sizes don't depend on waiting times, so that we may factorise the joint probability density $\phi(x,t)=\eta(t)w(x)$ in terms of the 2 marginal densities $\eta$ and $w$. A common choice for the initial condition would be to assume the agent to be situated at the origin $x=0$ at time $t=0$. In other words, let $p_0(x)=\delta(x)$. 

Once again, we will have to consider a suitable---diffusive---limit. 
%{\color{red}The limits $s\rightarrow0$ and $k \rightarrow0$ considered earlier can be recovered as a special case of the general description.} 
Let us first fix the time $t$ and let $N(t)$ be given by \eqref{Nt}. Scaling both the observed waiting times and jump lengths by multiplication with some small---but positive---parameters $r$ and $h$ respectively, we will derive the governing macroscopic (Fourier-Laplace) equation. Keep in mind that $r$ and $h$ will be related through some \textit{scaling relation} that will be made precise later on. 

Specifically, let us introduce the newly rescaled quantities, \cite{SubdiffFD}
% \begin{equation}
%     t_{N(t)}(r) = r\sum_{i=1}^{N(t)} \tau_i
% \end{equation}\label{tN}
% and
% \begin{equation}
%     x_{N(t)}(h) = h\sum_{i=1}^{N(t)} x_i.
% \end{equation}\label{xN}
\begin{equation}
    t_{N(t)}(r) = r\sum_{i=1}^{N(t)} \tau_i,\qquad x_{N(t)}(h) = h\sum_{i=1}^{N(t)} x_i.
\end{equation}\label{tNxN}
One of the basic properties of probability density functions regarding a change of variables allows us to introduce the PDFs $\eta_{r}$ and $w_h$ of the rescaled waiting times $(r\tau_i)_i$ and jump lengths $(hx_i)_i$ as
% \nkschange{
% \begin{align}
%     \eta_{r}(\tau)&=\frac{\eta(\tau/r)}{r}, \hspace*{0.2cm} (\tau>0) \hspace*{0.2cm} \text{ and }\label{phir}\\
%     w_h(x)&=\frac{w(x/h)}{h}.\label{wh}
% \end{align}
% }
\begin{equation}
    \eta_{r}(\tau)=\frac{\eta(\tau/r)}{r}, \  (\tau>0),\qquad
    w_h(x)=\frac{w(x/h)}{h}.
\end{equation}\label{phir-wh}
As for the Fourier and Laplace transforms, one simplifies 
\begin{align}
    \mathcal{L}[\eta_{r}](s)&=\frac{1}{r}\mathcal{L}[\eta(\tau/r)](s)%\nonumber\\
    =\frac{1}{r}r\mathcal{L}[\eta](rs)%\nonumber\\
    =\mathcal{L}[\eta](rs) \label{Lphir}
\end{align}
and
\begin{align}
    \mathcal{F}[w_h(x)]&=\frac{1}{h}\mathcal{F}[w(x/h)]%\nonumber \\
    =\frac{1}{h}|h|\mathcal{F}[w](hk)%\nonumber\\
    =\mathcal{F}[w](hk).\label{Fwh}
\end{align}
Altogether, the rescaled version of the master equation---in the Fourier-Laplace domain that is to say---is given by the modified (compared with \eqref{MontrollWeiss}), but still uncorrelated version of the Montroll-Weiss equation 
\begin{align}
    \mathcal{L}[\mathcal{F}[p](k,\cdot)](s)
    &=\frac{\mathcal{F}[p_0](k)}{s}\left\{\frac{1-\mathcal{L}[\eta_{r}](s)}{1-\mathcal{F}[w_h](k)\mathcal{L}[\eta_{r}](s)}\right\}=\frac{\mathcal{F}[p_0](k)}{s}\left\{\frac{1-\mathcal{L}[\eta](rs)}{1-\mathcal{F}[w](hk)\mathcal{L}[\eta](rs)}\right\}.\label{MontrollWeissscaled}
\end{align}
Note that \eqref{MontrollWeissscaled} is the Montroll-Weiss equation corresponding to a particle jumping a distance $hx_i$ after a waiting time $r\tau_i$. 

Recall that we have restricted ourselves to a specific class of waiting time- and jump densities which behave asymptotically according to power-law forms. We assume that $\eta$ satisfies the asymptotic relation
\begin{align}
    \eta(\tau)&\sim \tau^{-1-\beta}\label{asphi}
\end{align}
(with $\beta\in(0,1)$) while $w(x)$ takes on the specific form
\begin{align}
    w(x)&=[b+\epsilon(|x|)]|x|^{-(\alpha+1)}\label{asw}
\end{align}
with $\epsilon(|x|)\to 0$ as $|x|\rightarrow\infty$, for some coefficient $b>0$ and with exponent $\alpha\in (0,2)$, \cite{SubdiffFD}. 

Returning to the  Laplace and Fourier transforms, these can always be shown to be of the form
\begin{align}
    \mathcal{L}[\eta_{r}](s)&= \mathcal{L}[\eta](rs) %\nonumber\\
    = 1-c_1(rs)^{\beta}+\mathcal{O}(r^{\beta}) \hspace*{0.2cm} \text{ as $r\to 0$}\label{asLphi}
\end{align}
on the one hand and
\begin{align}
    \mathcal{F}[w_h](k)&=\mathcal{F}[w](hk)%\nonumber\\
    =1-c_2(h|k|)^{\alpha} +\mathcal{O}(h^{\alpha}) \hspace*{0.2cm} \text{ as $h\to 0$}\label{asFw}
\end{align}
on the other hand, \cite{SubdiffFD}. Note that $c_1$ and $c_2$ here represent certain constants. The details in the derivation of \eqref{asLphi} and \eqref{asFw} are more complicated than one might expect. The interested reader is referred to the final paragraph on page 5 (continued on page 6) in \cite{SubdiffFD}. Scalas et al. provide further guidance on the corollaries necessary to prove these general expressions \cref{asLphi}, \eqref{asFw} for $\mathcal{L}[\eta_{r}](s)$ and $\mathcal{F}[w_h](k)$, solely based on the asymptotic forms \eqref{asphi} and \eqref{asw}. 

Using \eqref{asLphi} and \eqref{asFw} in \eqref{MontrollWeissscaled} results in 
\begin{align}
    \mathcal{L}[\mathcal{F}[p](k,\cdot)](s)
    &=\frac{\mathcal{F}[p_0](k)}{s}\left\{\frac{1-[1-c_1(rs)^{\beta}+\mathcal{O}(r^{\beta})]}{1-[1-c_2(h|k|)^{\alpha} +\mathcal{O}(h^{\alpha})][1-c_1(rs)^{\beta}+\mathcal{O}(r^{\beta})]}\right\}\\
    &\sim \frac{\mathcal{F}[p_0](k)}{s}\left\{\frac{c_1(rs)^{\beta}}{c_2(h|k|)^{\alpha} +c_1(rs)^{\beta}}\right\}
    \label{MontrollWeissscaled2}
\end{align}
when looking at the leading order behaviour under the limits $r,h\to 0$. We will, however, not let $r$ and $h$ tend to zero in any arbitrary (independent) fashion. On the contrary, we will suppose that the limits $r,h\to 0$ are taken according to the scaling relation, \cite{SubdiffFD}
    \begin{equation}
        c_1r^{\beta}=c_2h^{\alpha}.\label{scaling}
    \end{equation}

The Hurst exponent $H=\frac{\beta}{\alpha}$ quantifies the (anomalous in general, for $\beta\neq 1$ and $\alpha \neq 2$) scaling. Note that $H=\frac{1}{2}$ (which is the classical parabolic scaling) amounts to ordinary diffusion. Applying the above scaling relation \eqref{scaling} in \eqref{MontrollWeissscaled2} results in the following asymptotic behaviour of $\mathcal{L}[\mathcal{F}[p](k,\cdot)](s)$ as $r,h\rightarrow0$
\begin{align}
    \mathcal{L}[\mathcal{F}[p](k,\cdot)](s)
    &\sim \mathcal{F}[p_0](k)\left\{\frac{s^{\beta-1}}{|k|^{\alpha} +s^{\beta}}\right\}\label{MontrollWeissfinal}\\
    &\sim \mathcal{L}[\mathcal{F}[u](k,\cdot)](s).
    \label{MontrollWeissscaled3}
\end{align}
The approximate asymptotical quantity \eqref{MontrollWeissscaled3} is precisely the Fourier-Laplace transform of the \textit{Green function} $u$, \cite{SubdiffFD}.  

\begin{definition}{Green function, \cite{SubdiffFD}.}
    The Green function $u$ is the fundamental solution of the fractional-diffusion Cauchy problem
% \begin{align}
%     \begin{array}{ll}
%     \displaystyle
%         \pdv[\beta]{}{t}u(x,t)=\pdv[\alpha]{}{|x|}u(x,t) & \hspace*{0.4cm}\alpha\in(0,2],\hspace*{0.1cm} \beta \in (0,1]\\
%         u(x,0^+)=\delta(x) & \hspace*{0.4cm}x\in(-\infty,\infty),\hspace*{0.1cm} t>0,\label{Green}
%     \end{array}
% \end{align}
\begin{align}
    \begin{cases}
    \displaystyle
        \pdv[\beta]{}{t}u(x,t)=\pdv[\alpha]{}{|x|}u(x,t), & \alpha\in(0,2],\ \beta \in (0,1],\\
        u(x,0^+)=\delta(x), & x\in(-\infty,\infty),\ t>0,\label{Green}
    \end{cases}
\end{align}

where $\displaystyle\pdv[\beta]{}{t}$ represents a fractional Caputo derivative while $\displaystyle\pdv[\alpha]{}{|x|}$ is a Riesz derivative. 
\end{definition}
Note that, in general, $p_0(x)$ can take on any form (i.e. is not necessarily equal to $\delta(x)$), in this case $u$ would simply have to satisfy $u(x,0^+)=p_0(x)$ instead of $u(x,0^+)=\delta(x)$ in \eqref{Green}.

\subsubsection{Fractional Derivatives}
We provide here a more detailed exposition of the various fractional derivative definitions used in \eqref{Green}. Both time- and space fractional derivatives are integro-differential operators. 

The Caputo time-fractional derivative of order $\beta$ ($m-1<\beta\leq m$, $m\in \mathbb{N}$) can be defined through its representation
\begin{equation}\label{seriesCaputo}
    \mathcal{L}\left [ \pdv[\beta]{}{t}f(t);s\right ] = s^{\beta}\Tilde{f}(s) - \sum_{k=0}^{m-1} s^{\beta-1-k}f^{(k)}(0^+)
\end{equation}
in the Laplace domain. In the above series expansion, $m$ represents the smallest natural number bounding $\beta$ from above, let us denote this by $m = \lceil \beta \rceil$, using the ceiling function $\lceil \cdot \rceil$. The explicit definition is then readily given by 
\begin{align*}
    \pdv[\beta]{}{t}f(t) &=\left \{ \begin{array}{ll}
        \displaystyle\frac{1}{\Gamma(m-\beta)}\int_0^t\frac{f^{(m)}(\tau)}{(t-\tau)^{\beta+1-m}}d\tau, &  ~\ m-1<\beta<m, \\
         \displaystyle\dv[m]{}{t}f(t), & ~\ \beta = m
    \end{array}\right.
\end{align*}
in \cite{analytic}. In the more restrictive case $0<\beta<1$, the time-fractional derivative is defined through:
\begin{definition}{Caputo fractional derivative ($0<\beta<1$), \cite{FracTheory}.}
For a sufficiently well-behaved function\footnote{\label{note1}{ Understood to be as smooth as required for the derivatives and integrals involved in the definition to make sense}} $f(t)$ and $0<\beta<1$, the Caputo fractional derivative of order $\beta$ can be defined by
    \begin{align}
        \dv[\beta]{}{t}f(t)
        &=\frac{1}{\Gamma(1-\beta)} \int_0^t \frac{\dv[]{}{\tau}f(\tau)}{(t-\tau)^{\beta}}d\tau \label{Caputo}\\
        &= \frac{1}{\Gamma(1-\beta)} \dv[]{}{t}\int_0^t \frac{f(\tau)}{(t-\tau)^{\beta}}d\tau - \frac{f(0^+)}{\Gamma(1-\beta)t^{\beta}}.\label{Caputo0}
    \end{align}
\end{definition}
%\nkscomm{we should say something about the "well behaved}

With respect to space, the \textit{Riesz derivative} naturally emerges. It is frequently put to use in space-fractional quantum mechanics, where it plays a significant role.   
\begin{definition}{Riesz derivative, \cite{FracTheory}.}
    The Riesz derivative is commonly referred to as the `pseudo-differential operator with symbol $-|k|^{\alpha}$', seeing as it is defined through its Fourier transform
    \begin{align*}
        \mathcal{F}\left[\pdv[\alpha]{g}{|x|}\right](k)&=-|k|^{\alpha}\mathcal{F}[g](k)
        \iff \pdv[\alpha]{g}{|x|}(k)=-\mathcal{F}^{-1}\left\{|k|^{\alpha}\mathcal{F}[g](k)\right\}.
    \end{align*}
    For $0<\alpha<2$, there also exists a regularised explicit integral representation, \cite{analytic}, which is given by
    %this means that the (symmetric) Riesz derivative 
    \begin{equation*}
        \pdv[\alpha]{}{|x|}g(x)=\frac{\Gamma(1+\alpha)}{\pi}\sin(\frac{\alpha\pi}{2})\int_0^{\infty}\frac{g(x-z)-2g(x)+g(x+z)}{z^{1+\alpha}}dz.
    \end{equation*}
\end{definition}
Riesz derivatives belong to the wider class of Riesz-Feller space-fractional derivatives. They make up the symmetric subclass of Riesz-Feller derivatives. 

The Riesz-Feller space-fractional derivative ${}_xD_{\theta}^{\alpha}$ of order $\alpha$, with skewness parameter $\theta$, is defined by means of its representation in the Fourier domain. Namely, for a sufficiently well-behaved function\footnoteref{note1} $g$, one has
\begin{equation}\label{fourrieszfeller}
    \mathcal{F}\left[{}_xD_{\theta}^{\alpha}g\right](k) = -\psi_{\alpha}^{\theta}(k)\mathcal{F}[g](k)
\end{equation}
The relation \eqref{fourrieszfeller} highlights in particular that ${}_xD_{\theta}^{\alpha}$ is a pseudo-differential operator with symbol $-\psi_{\alpha}^{\theta}(k)$. This symbol equals the logarithm of the characteristic function of a general L{\'e}vy strictly stable probability density, which has index of stability $\alpha$ and asymmetry parameter $\theta$, \cite{analytic}. As such,
\begin{equation}\label{psi}
    \psi_{\alpha}^{\theta}(k) = |k|^{\alpha}e^{i(\text{sign} k )\theta\pi/2}
\end{equation}
and, moreover, we continue to restrict the allowed region for the parameters $\alpha$ and $\theta$ to
\begin{equation}\label{dompars}
    0<\alpha\leq 2, ~\ |\theta|\leq \min\{\alpha,2-\alpha\}.
\end{equation}
In case $\theta=0$, the Riesz-Feller derivative becomes the ordinary Riesz space-fractional derivative. The space-fractional derivative then represents a symmetric operator with respect to the variable $x$. One denotes
\begin{equation}\label{riesz}
    {}_xD_{0}^{\alpha} = \dv[\alpha]{}{|x|}.
\end{equation}

Yet another representation of the Riesz derivative exists, in terms of left- and right Riemann-Liouville derivatives. 
Indeed, the Riesz fractional derivative can be constructed as the specific weighted sum of these two types of Riemann-Liouville derivatives
\begin{align}\label{defRiesz}
    \pdv[\alpha]{}{|x|}u(x,t) &= -c_{\alpha} \left( { }_{\text{LRL}}D_{a,x}^{\alpha} + { }_{\text{RRL}} D_{x,b}^{\alpha}\right) u(x,t),
\end{align}
with linear coefficient $c_{\alpha} = \frac{1}{2\cos(\frac{\alpha \pi}{2})}$ that depends on the order of the Riesz derivative, \cite{XieFang}. This linear combination of operators explicitly showcases the two-sided nature of the derivative. 

The integro-differential representation of left- and right Riemann-Liouville derivatives (taken from the fractional derivative definitions, for $\alpha>0$, on p. 1090 in \cite{BCtwosided}) is given by 
\begin{align*}
\begin{array}{llcl}
    \text{Left-side Riemann-Liouville derivative} & { }_{\text{LRL}}D_{a,x}^{\alpha} u(x,t) & = & \frac{1}{\Gamma(n-\alpha)} \dv[n]{}{x} \int_{a}^x \frac{u(y,t)}{(x - y)^{\alpha - n + 1}} dy,
    \\
    \text{Right-side Riemann-Liouville derivative} & { }_{\text{RRL}} D_{x,b}^{\alpha} u(x,t) & = & \frac{(-1)^n}{\Gamma(n-\alpha)} \dv[n]{}{x} \int_{x}^b \frac{u(y,t)}{(y - x)^{\alpha - n + 1}} dy. 
\end{array}
\end{align*}
Note that the natural number $n$ involved in the above definitions depends on $\alpha$, by means of $n = \lfloor \alpha \rfloor + 1$. $\Gamma(\cdot)$ is Euler's well-known Gamma function. 
Combining both integrals into one-and-the-same, 
\begin{align*}
    \pdv[\alpha]{}{x}u(x,t) &= \frac{c_{\alpha}}{\Gamma(n-\alpha)}\pdv[n]{}{x}\int_a^b \frac{u(y,t)}{|x-y|^{\alpha -1}}dy
\end{align*}
defines the Riesz derivative in an alternative way. 
We will further adopt the assumption that $1<\alpha<2$ (in which case the Riesz derivative exists), which makes that $n=2$ in the above. If $\alpha =2$, both left-and right Riemann-Liouville derivatives reduce to the ordinary derivative, by convention, \cite{BCtwosided}. For a comprehensive introduction to fractional derivatives, along with their properties, see also \cite{Podlubny1999}.
%{\color{blue}(ADD BOOKS OFFICE)}

\subsubsection{Green Function}
After defining the various types of fractional derivatives, we can check that $p$ and $u$ are indeed asymptotically equivalent in the Fourier-Laplace domain as \eqref{MontrollWeissscaled3} hypothesises. On the one hand, by definition,
\begin{align*}
    \mathcal{F}\left [ \pdv[\alpha]{}{|x|}u(x,t)\right](k) = -|k|^{\alpha} \mathcal{F}[u(\cdot,t)](k).
\end{align*}
Hence, also 
\begin{align}
    \mathcal{L}\left [\mathcal{F}\left [ \pdv[\alpha]{}{|x|}u(x,t)\right]\right](k,s) = -|k|^{\alpha} \mathcal{L}[\mathcal{F}[u(\cdot,t)](k)](s)\label{90}
\end{align}
on the Fourier-Laplace transform level.
On the other hand, by definition of the Caputo fractional derivative \eqref{Caputo0}, together with some basic identities regarding Laplace transforms (notably, the formula for the Laplace transform of a derivative), one finds
\begin{align}
    \mathcal{L}\left[\mathcal{F}\left[\pdv[\beta]{}{t}u(x,\cdot)\right]\right](k,s)
    &=\mathcal{L}\left[\pdv[\beta]{}{t}\mathcal{F}\left[u(x,\cdot)\right]\right](k,s)
    \nonumber\\
    &=\mathcal{L}\left[\frac{1}{\Gamma(1-\beta)} \dv[]{}{t}\int_0^t \frac{\mathcal{F}\left[u(\cdot,\tau)\right](k)}{(t-\tau)^{\beta}}d\tau - \frac{\mathcal{F}\left[u(\cdot,0^+)\right](k)}{\Gamma(1-\beta)t^{\beta}}\right](s)\nonumber\\
    &=\frac{1}{\Gamma(1-\beta)} s\mathcal{L}\left[\int_0^t \frac{\mathcal{F}\left[u(\cdot,\tau)\right](k)}{(t-\tau)^{\beta}}d\tau\right](s) - \mathcal{L}\left[\frac{\mathcal{F}\left[u(\cdot,0^+)\right](k)}{\Gamma(1-\beta)t^{\beta}}\right](s)\nonumber\\
    &=\frac{1}{\Gamma(1-\beta)} s\mathcal{L}\left[ \mathcal{F}\left[u(\cdot,t)\right](k)*_t(t)^{-\beta}\right](s) - \frac{\mathcal{F}\left[u(\cdot,0^+)\right](k)}{\Gamma(1-\beta)}\mathcal{L}\left[t^{-\beta}\right](s)\nonumber\\
     &=\frac{1}{\Gamma(1-\beta)} s\mathcal{L}\left[ \mathcal{F}\left[u(\cdot,t)\right](k)\right](s)\mathcal{L}\left[t^{-\beta}\right](s) - \frac{\mathcal{F}\left[u(\cdot,0^+)\right](k)}{\Gamma(1-\beta)}\mathcal{L}\left[t^{-\beta}\right](s)\nonumber\\
     &=\left\{ s\mathcal{L}\left[ \mathcal{F}\left[u(\cdot,t)\right](k)\right](s) - \mathcal{F}\left[u(\cdot,0^+)\right](k)\right\}\frac{1}{\Gamma(1-\beta)}\mathcal{L}\left[t^{-\beta}\right](s)\nonumber\\
     &=\left\{ s\mathcal{L}\left[ \mathcal{F}\left[u(\cdot,t)\right](k)\right](s) - \mathcal{F}\left[u(\cdot,0^+)\right](k)\right\}s^{\beta -1}.\label{98}
\end{align}
%Note that we have denoted $*_t$ to make clear that the convolution is taken with respect to the time variable. 

The function $u$ satisfies \eqref{Green} if and only if this equivalence holds at the level of Fourier-Laplace transforms. As a consequence of \eqref{90} and \eqref{98}, \eqref{Green} implies
\begin{align*}
    \left\{ s\mathcal{L}\left[ \mathcal{F}\left[u(\cdot,t)\right](k)\right](s) - \mathcal{F}\left[u(\cdot,0^+)\right](k)\right\}s^{\beta -1} =  -|k|^{\alpha} \mathcal{L}[\mathcal{F}[u(\cdot,t)](k)](s).
\end{align*}
Or, equivalently, 
\begin{align*}
    \left(s^{\beta} +|k|^{\alpha} \right)\mathcal{L}\left[ \mathcal{F}\left[u(\cdot,t)\right](k)\right](s) &= \mathcal{F}\left[u(\cdot,0^+)\right](k)s^{\beta -1} \\
    \iff \mathcal{L}\left[ \mathcal{F}\left[u(\cdot,t)\right](k)\right](s) &= \mathcal{F}\left[u(\cdot,0^+)\right](k)\frac{s^{\beta -1} }{s^{\beta} +|k|^{\alpha}}
\end{align*}
which is precisely the exact equality resulting from \eqref{MontrollWeissfinal}, for $p=u$. Adopting the specific initial condition $u(x,0^+)=\delta(x)$, we have justified the postulated Fourier-Laplace asymptotic equivalence \eqref{MontrollWeissscaled3} between the probability density function $p$ and the Green function as defined by \eqref{Green}.

As a result, $p$ and $u$ will be weakly asymptotically equivalent in the space-time domain. Indeed, note that we have adopted the non-unitary, angular frequency definition for the Fourier transform.  This means that
\begin{align}
    \mathcal{L}\left[\mathcal{F}\left[p\right]\right](k,s)&=\int_0^{\infty} \mathcal{F}[p(\cdot,t)](k)e^{-st}dt
    =\int_0^{\infty} \int_{\R} p(x,t)e^{ikx}e^{-st}~dx dt\label{LF}
\end{align}
whereas the characteristic function of the \emph{random vector} $(x,t)$ is given by
\begin{align}
    \varphi_{p,(x,t)}(k,s) &= \int_0^{\infty} \int_{\R} p(x,t)e^{ikx}e^{ist}~dxdt. \label{char}
\end{align}
Upon direct comparison, it can be seen that 
\begin{equation}
    \mathcal{L}\left[\mathcal{F}\left[p\right]\right](k,-is) = \varphi_{p,(x,t)}(k,s).
\end{equation}
Upon taking the limits $r,h\to 0$ (according to the scaling relation \eqref{scaling}), \eqref{MontrollWeissscaled3} tells us precisely that $\mathcal{L}\left[\mathcal{F}\left[p\right]\right](k,s)\to \mathcal{L}\left[\mathcal{F}\left[u\right]\right](k,s)$ pointwise, for any real $k$ and complex $s$. Hence, also $\varphi_{p,(x,t)}(k,s)=\mathcal{L}\left[\mathcal{F}\left[p\right]\right](k,-is)\to \mathcal{L}\left[\mathcal{F}\left[u\right]\right](k,-is)=\varphi_{u,(x,t)}(k,s)$ pointwise, on the level of characteristic functions. L{\'e}vy's Continuity Theorem, for sequences of characteristic functions, then tells us precisely that $p$ weakly converges to $u$. In this sense, we can say that the master equation \eqref{massprob} is asymptotically equivalent (in a weak sense) to the fractional diffusion equation \eqref{Green} in the space-time domain, per application of L{\'e}vy's Continuity Theorem.

\subsection{Mesoscopic Density Equations}
%{\color{red} Does it make sense to include a multidimensional discussion here? Or should we simplify to the 1D case?}\nkscomm{I think it's ok to leave it as is}
In the last section, we were able to derive a fractional diffusion equation for the probability density function of a single agent, in the long-time large-scale limit. However, since we are interested in the collective migratory behaviour of multiple agents, we are hoping to find similar results for the density $\rho$. We here opt for a mesoscopic instead of an individualistic approach, where one introduces mean-field equations for the particle density. The mesoscopic equations themselves will, however, involve a detailed description of the particle movement at microscopic level. Instead of restricting ourselves to the one-dimensional scenario, we will at this point look at the more general $d$-dimensional setup (where $d=1,2,3$ are the scenarios of biological relevance), since the underlying mathematics are mostly the same. Every random step undertaken by a single particle constituting to the density $\rho$ is still characterised by both the waiting time and jump length, as distributed according to the joint probability function $\phi(\Vec{y},t)$, with $\Vec{y}\in \R ^d$ and $t\in [0,\infty)$. 

Assuming the particles follow a CTRW (individually), one derives the mean-field equations mentioned in \cite{reaction-transport}
\begin{align}
    \underbrace{\rho(\Vec{x},t)}_{\substack{\text{particle density at} \\ \text{position $\Vec{x}$ at time $t$}}}&=\underbrace{\rho(\Vec{x}, 0)}_{\substack{\text{initial particle density } \\ \text{at the position $\Vec{x}$ at time $0$ }}}\underbrace{\Psi(t)}_{\substack{\text{ probability of no jump occurring}\\ \text{ after a waiting time less than $t$}}} \nonumber\\
    & \hspace{0.2cm} + \int_0^t \underbrace{j(\Vec{x}, t-\tau)}_{\substack{\text{density of particles }\\ \text{present at the position }\\ \text{$\Vec{x}$ at time $t-\tau$ }}}\underbrace{\Psi(\tau)}_{\substack{\text{ probability of the particles not}\\ \text{ leaving their position for a }\\ \text{ duration of length at least $\tau$}}}d\tau,\label{rho}
\end{align}
with
\begin{align}
    j(\Vec{x},t)&=\underbrace{F(\rho(\Vec{x}, t))}_{\substack{\text{density of particles }\\ \text{being produced at } \\ \text{ position $\Vec{x}$ and time $t$ }}} + \underbrace{\int_{\R ^3} \rho(\Vec{x}-\Vec{z}, 0)\phi(\vec{z},t)~d\Vec{z}  + \int_0^t\int_{\R ^3} j(\Vec{x}-\Vec{z},t-\tau)\phi(\Vec{z},\tau)~d\Vec{z}d\tau}_{\substack{\text{ density of particles that arrive at the position $\Vec{x}$ exactly at time $t$,} \\ \text{ from other points $\Vec{x}-\Vec{z}$, by means of transport}}}. \label{jreact}
\end{align}
Note that we have now included two processes according to which particles (dis)appear at position $\Vec{x}$. The first is due to particle transport, while the second can be attributed to creation/destruction phenomena. The latter is captured by the density-dependent rate $F(\rho(\Vec{x}, t))$ of density being produced at the point $\vec x$ at time $t$. In the density setting, $j(\vec x,t)$ can be seen as ``the new density being added to the existing'', \cite{reaction-transport}.

Comparing \eqref{massprob} and \eqref{jprob} to \eqref{rho} and \eqref{jreact}, the density equations take the same form as the prior probability density ones. In fact, all one has to do is exchange $p$ for $\rho$ in \eqref{massprob} and \eqref{jprob}, and add $F(\rho(\Vec{x},t))$ to the resulting expression for $j(\Vec{x},t)$. Note also that $\Psi(t)$ is still defined as the survival probability \eqref{waitingtime}. Once again, we are interested in the long-time large-scale limit. This time, it will lead us to what will be identified as a fractional reaction-transport equation. However, we first need to find a density-specific reaction-diffusion analogue to the Montroll-Weiss equation \eqref{MontrollWeiss}. The ultimate goal is to find a more comprehensible expression for $\rho(x,t)$ as defined by \eqref{rho}. As before, we will do so by considering Fourier-Laplace transforms. 

Taking the Fourier transform $\mathcal{F}$---with respect to the space vector variable $\Vec{x}$---of \eqref{rho} leads to (in terms of the new $d$-dimensional variable $\Vec{k}$ of the transform)
\begin{align}
    \mathcal{F}[\rho](\Vec{k},t) &= \mathcal{F}[\rho_0](\Vec{k})\Psi(t)+\int_0^t\mathcal{F}[j](\Vec{k},t-\tau)\Psi(\tau)d\tau=\mathcal{F}[\rho_0](\Vec{k})\Psi(t)+\mathcal{F}[j](\Vec{k},\cdot)*_t\Psi(\cdot), \label{fourierrho}
\end{align}
which is completely analogous to what we had previously found in \eqref{fourierp}. 

Consecutively taking the Laplace transform of \eqref{fourierrho} (moving from the time variable $t$ to the Laplace variable $s$, whilst keeping in mind multiplicativity of the aforementioned transform with respect to convolutions) then brings us to 
\begin{align}
    \mathcal{L}[\mathcal{F}[\rho](\Vec{k},\cdot)](s)
    &=\mathcal{F}[\rho_0](\Vec{k})\mathcal{L}[\Psi](s)+\mathcal{L}[\mathcal{F}[j](\Vec{k},\cdot)*_t\Psi(\cdot)](s) \nonumber\\
    &=\mathcal{F}[\rho_0](\Vec{k})\mathcal{L}[\Psi](s)+\mathcal{L}[\mathcal{F}[j](\Vec{k},\cdot)](s)\mathcal{L}[\Psi(\cdot)](s).\label{FLrho1}
\end{align}
The only difference upon including creation effects in the density mass equation appears in the formula for $\mathcal{L}[\mathcal{F}[j](\Vec{k},\cdot)](s)$, where we just have to modify the multidimensional version of \eqref{fourierjtotal} for $\rho$, by adding $\mathcal{F}[F(\rho(\cdot,t))](\Vec{k})$ to the right-hand side. Hence, 
\begin{align}
    \mathcal{F}[j](\Vec{k},t) &= \mathcal{F}[\phi](\Vec{k},t)\mathcal{F}[\rho_0](\Vec{k}) + \mathcal{F}[j](\Vec{k},\cdot)*_t\mathcal{\phi}(\Vec{k},\cdot) + \mathcal{F}[F(\rho(\cdot,t))](\Vec{k}).\label{fourierjtotalrho}
\end{align}
This means that, in terms of the Fourier-Laplace transform of $j$, we now find
\begin{align}
    \mathcal{L}[\mathcal{F}[j]](\Vec{k},s) &= \mathcal{L}[\mathcal{F}[\phi]](\Vec{k},s)\mathcal{F}[\rho_0](\Vec{k}) + \mathcal{L}[\mathcal{F}[j]](\Vec{k},s)\mathcal{L}[\mathcal{F}[\phi]](\Vec{k},s) + \mathcal{L}[\mathcal{F}[F]](\Vec{k},s).
\end{align}
Therefore, altogether
\begin{align}
    \mathcal{L}[\mathcal{F}[j]](\Vec{k},s) &=\frac{\mathcal{L}[\mathcal{F}[\phi]](\Vec{k},s)\mathcal{F}[\rho_0](\Vec{k}) + \mathcal{L}[\mathcal{F}[F]](\Vec{k},s)}{1-\mathcal{L}[\mathcal{F}[\phi]](\Vec{k},s)} \label{laplacefourierjrho}
\end{align}
in terms of the---presumably known---initial density $\rho_0$ and CTRW's characteristic PDF $\phi$. 
Note that $\mathcal{L}[\Psi](s)$ remains unchanged and equal to \eqref{laplacepsifinal} if we continue to assume independence also in this multi-dimensional scenario. That is to say, jump lengths and waiting times are assumed independent of one another, allowing us to factorise the joint PDF $\phi(\Vec{y},t)=w(\Vec{y})\eta(t)$ in terms of the marginal waiting time probability density function $\eta(t)$ along with the jump length probability density function $w(\Vec{y})$.
Inserting the expressions \eqref{laplacefourierjrho} and \eqref{laplacepsifinal} into \eqref{FLrho1}, one finds the reaction-diffusion \textit{Montroll-Weiss} analogue satisfied by the density $\rho$. Namely,
\begin{align}
    \mathcal{L}[\mathcal{F}[\rho](\Vec{k},\cdot)](s)
    &=\left\{\mathcal{F}[\rho_0](\Vec{k})+\mathcal{L}[\mathcal{F}[j](\Vec{k},\cdot)](s)\right\}\mathcal{L}[\Psi(\cdot)](s)\nonumber\\
    &=\left\{\mathcal{F}[\rho_0](\Vec{k})+\frac{\mathcal{L}[\mathcal{F}[\phi]](\Vec{k},s)\mathcal{F}[\rho_0](\Vec{k})+\mathcal{L}[\mathcal{F}[F]](\Vec{k},s)}{1-\mathcal{L}[\mathcal{F}[\phi]](\Vec{k},s)}\right\}\frac{1-\mathcal{L}[\eta](s)}{s} \nonumber\\
    &=\left\{\frac{\mathcal{F}[\rho_0](\Vec{k})+\mathcal{L}[\mathcal{F}[F]](\Vec{k},s)}{1-\mathcal{L}[\mathcal{F}[\phi]](\Vec{k},s)}\right\}\frac{1-\mathcal{L}[\eta](s)}{s}.\label{MontrollWeissrho}
\end{align}
\subsubsection{(Sub/Super)diffusive Limit}
As in \eqref{tNxN}, space and time should be rescaled by factors $h$ respectively $r$ (the only difference being that the $x_i$ are now replaced by their vector analogues $\Vec{x}_i$). On the level of Fourier and Laplace transforms, \eqref{Lphir} and the vector analogue of \eqref{Fwh} together lead to the (re)scaled Montroll-Weiss equation
\begin{align}
    \mathcal{L}[\mathcal{F}[\rho](\Vec{k},\cdot)](s)
    &=\left\{\frac{\mathcal{F}[\rho_0](\Vec{k})+\mathcal{L}[\mathcal{F}[F]](\Vec{k},s)}{1-\mathcal{F}[w](h\Vec{k})\mathcal{L}[\eta](rs)]}\right\}\frac{1-\mathcal{L}[\eta](rs)}{s}
\end{align}
for (the Fourier-Laplace transform of) the density $\rho$. 
Assume asymptotic heavy-tailed behaviour 
\begin{align*}
    w(\Vec{y})\sim \|\Vec{y} \|^{-(\alpha + d)} \quad \text{ as }\quad \|\Vec{y}\|\to \infty, \qquad \eta(\tau)\sim \tau^{-1-\beta} \quad \text{ as }\quad \tau \to \infty,
\end{align*}
where the $d$ in the negative exponent stems from the $d$-dimensionality of the underlying space, with $\alpha\in (0,2)$ and $\beta \in (0,1)$, \cite{SubdiffFD}.

\textbf{For $\mathbf{d=1}$,} the same approximations \eqref{asLphi} and \eqref{asFw} are retrieved. This results in the asymptotic behaviour
\begin{align}
    &\mathcal{L}[\mathcal{F}[\rho](k,\cdot)](s)\\
    &=\left\{\frac{\mathcal{F}[\rho_0](k)+\mathcal{L}[\mathcal{F}[F]](k,s)}{1-\left\{1-c_2(h|k|)^{\alpha} +\mathcal{O}(h^{\alpha})\right\}\left\{1-c_1(rs)^{\beta}+\mathcal{O}(r^{\beta})\right\} }\right\}\frac{1-\left\{1-c_1(rs)^{\beta}+\mathcal{O}(r^{\beta})\right\} }{s}\nonumber\\
    &\sim\left\{\frac{\mathcal{F}[\rho_0](k)+\mathcal{L}[\mathcal{F}[F]](k,s)}{c_2(h|k|)^{\alpha} + c_1(rs)^{\beta}}\right\}\frac{c_1(rs)^{\beta}}{s}
    \label{MontrollWeissscaledrho}
\end{align}
as $r,h\to 0$. 
%\nkschange{I have removed the numbering. Consider repeating this "inline" rather than in "display"}
We derive the limiting behaviour
\begin{align}
    \mathcal{L}[\mathcal{F}[\rho](k,\cdot)](s)
    &\sim\left\{\mathcal{F}[\rho_0](k)+\mathcal{L}[\mathcal{F}[F]](k,s)\right\}\frac{s^{\beta-1}}{|k|^{\alpha} + s^{\beta}}
    \label{MontrollWeissscaledrhofinal}
    \\ &\sim \mathcal{L}[\mathcal{F}[u](k,\cdot)](s)
\end{align}
as $r,h\rightarrow0$, according to \eqref{scaling} as before.
This time, the density $\rho$ can be seen to be weakly asymptotically equivalent to the function $u:\R\times(0,\infty)\rightarrow\R$ which satisfies the following, fractional, reaction-diffusion equation
    \begin{align}
        \pdv[\beta]{}{t}u(x,t)=\pdv[\alpha]{}{|x|}u(x,t) + \frac{1}{\Gamma(1-\beta)}\int_0^t \frac{F(x,\tau)}{(t-\tau)^{\beta}}d\tau, \quad \alpha\in(0,1)\cup(1,2),\ \beta \in (0,1), \label{fracrho}
    \end{align}
subject to some initial condition
\begin{align}
    u(x,0^+)=\rho(x,0^+). \label{rhoinitial}
\end{align}
Upon taking the Fourier-Laplace transform of the \textit{fractional reaction-diffusion equation} \eqref{fracrho}, one readily computes the solution $u$ of \eqref{fracrho} (subject to \eqref{rhoinitial}) to satisfy
\begin{equation}
    \mathcal{L}[\mathcal{F}[u](k,\cdot)](s) = \left\{\mathcal{F}[\rho_0](k)+\mathcal{L}[\mathcal{F}[F]](k,s)\right\}\frac{s^{\beta-1}}{|k|^{\alpha} + s^{\beta}}. 
\end{equation}
In fact, this can be justified by looking at the F-L transforms of the individual components that contribute to the fractional equation \eqref{fracrho}. The equations \eqref{90} and \eqref{98} are still valid, namely
\begin{align}
    \mathcal{L}\left [\mathcal{F}\left [ \pdv[\alpha]{}{|x|}u(x,t)\right]\right](k,s) &= -|k|^{\alpha} \mathcal{L}[\mathcal{F}[u(\cdot,t)](k)](s)
\end{align}
and
\begin{align}
    \mathcal{L}\left[\mathcal{F}\left[\pdv[\beta]{}{t}u(x,\cdot)\right]\right](k,s)
     &=\left\{ s\mathcal{L}\left[ \mathcal{F}\left[u(\cdot,t)\right](k)\right](s) - \mathcal{F}\left[u(\cdot,0^+)\right](k)\right\}s^{\beta -1}.
\end{align}
As for the `fractional reaction term', one turns to some basic rules of computation with Fourier- and/or Laplace transforms
\begin{align*}
    \mathcal{L}\left[\mathcal{F}\left[\frac{1}{\Gamma(1-\beta)}\int_0^t \frac{F(x,\tau)}{(t-\tau)^{\beta}}d\tau\right]\right](k,s) &= \mathcal{L}\left[\mathcal{F}\left[\frac{1}{\Gamma(1-\beta)}F(x,t)*_t(t^{-\beta})\right]\right](k,s)\\
    &= \mathcal{L}\left[\frac{1}{\Gamma(1-\beta)}\mathcal{F}\left[F(\cdot,t)\right](k)*_t(t^{-\beta})\right](s),
\end{align*}
where one recalls that the Fourier transform is taken with respect to the $x$ variable. Furthermore,
\begin{align*}
    \mathcal{L}\left[\frac{1}{\Gamma(1-\beta)}\mathcal{F}\left[F(\cdot,t)\right](k)*_t(t^{-\beta})\right](s) &= \frac{1}{\Gamma(1-\beta)}\mathcal{L}\left[\mathcal{F}\left[F(\cdot,t)\right](k)\right](s)\mathcal{L}\left[t^{-\beta}\right](s)
\end{align*}
due to the multiplicative property of Fourier transforms (with respect to the time variable $t$). Moreover, the Laplace transform of the power function $t^{-\beta}$---with $0\leq\beta <1$---equals
\begin{align*}
    \mathcal{L}\left[t^{-\beta}\right](s) &= s^{\beta-1} \Gamma(1-\beta),
\end{align*}
if $\Re{s}>0$.
%(where this function is only assumed to be defined for $t\geq 0$, with $t$ representing the temporal variable). 
Then, one finds
\begin{align*}
    \mathcal{L}\left[\frac{1}{\Gamma(1-\beta)}\mathcal{F}\left[F(\cdot,t)\right](k)*_t(t^{-\beta})\right](s) &= \mathcal{L}\left[\mathcal{F}\left[F(\cdot,t)\right](k)\right](s)s^{\beta-1}.
\end{align*}
Altogether, combining the individual F-L transforms, the Fourier-Laplace equation resulting from \eqref{fracrho} boils down to the following equivalences
\begin{align*}
% \begin{array}{lrcl}
    \left\{ s\mathcal{L}\left[ \mathcal{F}\left[u(\cdot,t)\right](k)\right](s) - \mathcal{F}\left[u(\cdot,0^+)\right](k)\right\}&s^{\beta -1} \\
     =& -|k|^{\alpha} \mathcal{L}[\mathcal{F}[u(\cdot,t)](k)](s) + \mathcal{L}\left[\mathcal{F}\left[F(\cdot,t)\right](k)\right](s)s^{\beta-1}
    \\
    \iff\mathcal{L}\left[ \mathcal{F}\left[u(\cdot,t)\right](k)\right](s) \left \{s^{\beta}+|k|^{\alpha}\right \} =& \{\mathcal{L}\left[\mathcal{F}[u(\cdot,t)](k)\}(s)+\mathcal{F}\left[u(\cdot,0^+)\right](k)\right]s^{\beta -1} \\
    \iff\mathcal{L}\left[ \mathcal{F}\left[u(\cdot,t)\right](k)\right](s)  =& \{\mathcal{L}\left[\mathcal{F}\left[u(\cdot,0^+)+\mathcal{F}[u(\cdot,t)](k)\}(s)\right](k)\right]\frac{s^{\beta -1}}{s^{\beta}+|k|^{\alpha}},
% \end{array}
\end{align*}
with the latter being \eqref{MontrollWeissscaledrhofinal} precisely. Hence, $\rho$ actually (weakly) converges to $u$ in the long-time, large-scale limit. We have thus deduced the macroscopic fractional reaction-transport equation \eqref{fracrho}, starting from the basic CTRW density mass conservation equations \eqref{rho} and \eqref{jreact}. 

%======================================================================
\section{A Fractional Stochastic Differential Equation}\label{sec:fractSDE}
Let us go back to the setup described in Section \ref{CTRWs}. Instead of characterising the resulting probability density function $p(x,t)$ of a particle being at position $x$ at time $t$ through its integral master equation, it is possible to come up with a more explicit solution. This solution takes the form of a series that runs over the number of jumps $n$ occurring up until time $t$. In the uncoupled case, which we are still mostly invested in, the probability $P(n,t)$ of $n$ jumps occurring up to time $t$ plays a key part in the summation term, \cite{SubdiffFD}. If we assume that $p_0(x) = \delta(x)$, we recover
\begin{equation}\label{seriesp}
    p(x,t) = \sum_{n=0}^{\infty} P(n,t)w_n(x), ~\ \text{with} ~\ w_n(x) = \underbrace{(w * w * \hdots * w)}_{n \text{ times}}(x)
\end{equation}
the probability density of the sum $\sum_{i=1}^n x_i$.

%\subsection{Exponential Waiting-Time \cite{SubdiffFD}}
%In the last section, we managed to derive a general series description \eqref{seriesp} for the probability density function of a jump process, subordinate to a renewal process. Here, we provide an explicit example.
In case the waiting time is given to be of exponential form  $\eta(t) = \mu e^{-\mu t}$ ($\mu>0$), $\Psi(t) = e^{-\mu t}$, $P(n, t)$ equals the Poisson distribution
\begin{equation}
    P(n, t) =  \frac{(\mu t)^{n}}{n!}e^{-\mu t}.
\end{equation}
If, on the other hand, the survival probability function follows power-law decay 
\begin{equation}
    \Psi(\tau) \approx \tau^{-\beta}, ~\ \text{ as } \tau \to \infty ~\ \text{($0<\beta<1$)},
\end{equation}
one retrieves
\begin{equation}
    \Psi(\tau) = E_{\beta}(-\tau^{\beta}) ~\ \text{as defined by} ~\ E_{\beta}(z) = \sum_{n=0}^{\infty} \frac{z^n}{\Gamma(\beta n +1)}.
\end{equation}
This complex power series acts as the desired generalisation, which includes anomalous relaxation. In this case,
\begin{align}
    P(n,t) = \frac{t^{\beta n}}{n!}\dv[n]{}{z}E(z)\bigg|_{z=-t^{\beta}}.
\end{align}
In particular, for $\beta=1$, $P(n,t)$ reduces to the Poisson distribution again, seeing as Mittag-Leffler- and exponential functions coincide in this case, namely
\begin{equation}
    E(-t) = e^{-t}. 
\end{equation}

\subsection{A Comparison between SDEs and FSDEs}
Adopting an agent-based modelling approach, we are in need of a generalised form for a stochastic differential equation (SDE), which incorporates fractional diffusion. Starting from the foundations, we know that a Wiener process underlies most ordinary SDEs modelling regular diffusion. 
%Here, Brownian motion gives rise to regular diffusion. 
More generally speaking, L{\'e}vy processes underlie what we call \textit{FSDEs}. Essentially being a continuous-time analogue of the random walk, the distribution between successive displacements of the same particle takes on a more general form than that of a Gaussian distribution. 
%The setting is still widely the same as is the case for Wiener processes, with them actually representing a specific subset of the more general L{\'e}vy processes. 
%Notably, a L{\'e}vy process is a stochastic process which captures the motion of a particle undergoing random successive displacements. Moreover, the displacements in pairwise disjoint time intervals are independent, and displacements in distinct time intervals of the same length follow identical distributions. 

In the Brownian motion, regular diffusion case that we investigated previously, the probability of $n$ jumps occurring in a time interval of length $t$ was characterised by a Poisson distribution 
\begin{equation*}
    P(n,t) = \frac{(\mu t)^n}{n!}e^{-\mu t},
\end{equation*}
which has expected number of jumps $\lambda = \mu t$ within this time frame of length $t$. The probability density 
\begin{equation*}
    p(x,t) = \sum_{n=0}^{\infty}\frac{(\mu t)^n}{n!}e^{-\mu t}w_n(x)
\end{equation*}
could be seen to satisfy the regular diffusion equation 
\begin{equation}\label{pdvregdiff}
    \pdv{p}{t} = D\pdv[2]{p}{x}, ~\ D = \frac{\sigma^2}{2 \overline{T}},
\end{equation}
in the macroscopic limit. Denoting $X_t^i$ to be the position of the $i$-th particle at time $t$, it can be seen that the distribution of $X_{t+\Delta t}^i - X_t^i$ (equal to the distribution of $X_{\Delta t}^i$) is given by the Gaussian density
\begin{equation*}
    p(x,\Delta t) = \frac{1}{\sqrt{4\pi D \Delta t}}e^{-\frac{x^2}{4D\Delta t}}.
\end{equation*}
This is a consequence of the equation \eqref{pdvregdiff} which has an analytical solution equal to
\begin{equation*}
    p(x, t) = \frac{1}{\sqrt{4\pi D t}}e^{-\frac{x^2}{4Dt}}.
\end{equation*}

As for fractional-diffusive behaviour, we have mentioned before that 
\begin{equation*}
    P(n,t) = \frac{t^{\beta n}}{n!}\dv[n]{}{z}E_{\beta}(z)\bigg|_{z=-t^{\beta}}
\end{equation*}
generalises the Poisson distribution. In the macroscopic limit, 
\begin{equation}
    p(x,t) = \sum_{n=0}^{\infty}\frac{t^{\beta n}}{n!}\dv[n]{}{z}E_{\beta}(z)\bigg|_{z=-t^{\beta}}w_n(x)
\end{equation}
tends to the fundamental solution---i.e. Green's function $u(x,t)$---of the fractional diffusive equation
\begin{equation}
    \pdv[\beta]{}{t}u(x,t) = \pdv[\alpha]{}{|x|}u(x,t) ~\ \text{ for } ~\ x \in \R, t\in (0,\infty) .
\end{equation}
Hence, the distribution of $X_{t+\Delta t}^i - X_t^i$ in this case 
coincides with $u(x,\Delta t)$, in the limit. The subsequent paragraph will give an explicit calculation of $u(x, \Delta t)$, focusing on the $\beta =1$ (and general $\alpha$) case. As it turns out, 
\begin{equation}\label{distrXti}
    X_{\Delta t}^i \sim (\Delta t)^{-\beta/\alpha}W_{\alpha, \beta}\left(\frac{X}{\Delta t^{\beta/\alpha}}\right),
~\ \text{where} ~\
    W_{\alpha, \beta}(u) = \frac{1}{2\pi}\int_{\R} e^{-iku}E_{\beta}(-|k|^{\alpha})dk
\end{equation}
is the inverse Fourier transform of a Mittag-Leffler function.

\subsection{The Analytical Solution for Space-Fractional Diffusion}
In its most general form, the Cauchy problem for the space-time fractional diffusion equation is given by
\begin{equation}\label{fraccauchy}
    {}_xD_{\theta}^{\alpha}u(x,t) = \pdv[\beta]{}{t}u(x,t) ~\ \text{ for } ~\ x \in \R  \text{ and } t \in \R ^+,
\end{equation}
where $u(x,t)$ represents a (real) field variable depending on space ($x$) and time ($t$). 
It mimics the standard diffusion equation, where the second order space derivative is now replaced with a Riesz-Feller derivative ${}_xD_{\theta}^{\alpha}$ of order $\alpha \in (0,2]$ on the left-hand side of the equation \eqref{fraccauchy}. In addition, the Riesz-Feller derivative has a skewness parameter, represented by $\theta$ (where $|\theta|\leq \min\{\alpha, 2-\alpha\}$). However, we are only interested in the symmetric case here, which coincides with $\theta = 0$. On the right-hand side of \eqref{fraccauchy}, one notices that the first-order time derivative is replaced by the Caputo derivative $\pdv[\beta]{}{t}$ of general order $\beta\in(0,2]$, \cite{analytic}.

\paragraph{Fundamental solution}
The corresponding fundamental solution of the Cauchy problem \eqref{fraccauchy}, 
\begin{align}
   \text{ augmented with the initial condition }~\ u(x,0) = \delta(x), ~\ x \in \R ,\label{initialcond}
\\
     \text{ in addition to the boundary conditions } ~\ u(\pm\infty,t)=0, ~\ t>0, \label{BC}
\end{align}
is given by 
\begin{equation}
        u(x,t) = G_{\alpha, \beta}^{\theta}(x,t) = t^{-H}K^{\theta}_{\alpha,\beta}(x/t^{H}),\label{fundfrac}
\end{equation}
where we set $H = \beta/\alpha$ equal to the Hurst exponent.
One calls $K^{\theta}_{\alpha,\beta}$ the \textit{reduced Green function}. For more general, sufficiently well-behaved initial conditions $\varphi_0\in L^c(\R )$
\begin{equation}\label{init2}
      u(x,0) = \varphi_0(x), ~\ x \in \R ,
\end{equation}
the solution $u_{\alpha,\beta}^{\theta}$ of the fractional diffusion equation \eqref{fraccauchy}, subject to \eqref{BC} and \eqref{init2}, can be represented as the convolution of the initial condition $\varphi_0$ with the fundamental solution, the Green function $G_{\alpha,\beta}^{\theta}(x,t)$. One may write
\begin{equation}\label{gensol}
    u_{\alpha,\beta}^{\theta}(x,t) = \int_{-\infty}^{\infty} G_{\alpha, \beta}^{\theta}(\xi, t)\varphi_0(x-\xi)d\xi.
\end{equation}
For a more detailed explanation, the reader is referred to \cite{analytic}. 

It should be noted that, for $k\in \R $, $t\geq 0$, a representation of the fundamental solution can be given by
\begin{equation}\label{Greenexpr}
    G_{\alpha,\beta}^{\theta}(x,t)=\mathcal{F}^{-1}\left\{E_{\beta}\left[-\psi_{\alpha}^{\theta}(k)t^{\beta}\right]\right\}(x).
\end{equation}
Here, $E_{\beta}$ denotes the entire and transcendental Mittag-Leffler function of order $\beta$. This function is defined in the whole complex plane, for all non-negative values of the parameter $\beta$. Moreover, it can be expressed explicitly through the power series
\begin{equation}\label{MLfct}
    E_{\beta}(z) = \sum_{n=0}^{\infty} \frac{z^n}{\Gamma(\beta n +1)}, ~\ \beta>0, ~\ z \in \mathbb{C}.
\end{equation}
In order for the Fourier inverse of the function $E_{\beta}\left[-\psi_{\alpha}^{\theta}(k)t^{\beta}\right]$ in \eqref{Greenexpr} to exist, one further requires $|\theta|\leq2-\beta$, \cite{analytic}.

\paragraph{$\bm{\beta=1}$}
We here focus on the particular case of \textit{space-fractional diffusion} in which $\beta$ reduces to the regular-case value of one, while we remain free to vary $\alpha$ between $0<\alpha\leq 2$. This specific case $0<\alpha\leq 2, ~\ \beta =1$ of space-fractional diffusion allows us to interpret the fundamental solution as a spatial probability density function which evolves in time, \cite{analytic}.
If we set $\alpha = 2$, standard diffusion is recovered. The Green function is now given by the Gaussian probability density function
\begin{equation}
    G_{2,0}^{0}=\frac{1}{2\sqrt{\pi t}}e^{\frac{-x^2}{4t}}.\label{Gausspdf}
\end{equation}
The corresponding variance is evidently equal to
\begin{equation}
    \sigma^2 = 2t \label{vargauss}
\end{equation}
and is thus proportional to the first power of time (see also Einstein's diffusion law), \cite{analytic}. 

In the more general case $0<\alpha\leq 2$, the fact that $\beta=1$ implies that the Mittag-Leffler function $E_{1}$ reduces to an exponential function (according to \eqref{MLfct}, along with the Taylor expansion representation of the exponential function). Hence, $G_{\alpha,1}^{\theta}$ has to fulfil
\begin{equation}\label{FG}
    \mathcal{F}[G_{\alpha,1}^{\theta}](k,t) = e^{-t\psi_{\alpha}^{\theta}(k)}.
\end{equation}
If one denotes by $L_{\alpha}^{\theta}(x)$ the L{\'e}vy strictly stable density with parameters $\alpha$ and $\theta$, note that the characteristic function of a L{\'e}vy random variable is given by
\begin{equation}\label{FL}
    \mathcal{F}[L_{\alpha}^{\theta}](k) = e^{-\psi_{\alpha}^{\theta}(k)}.
\end{equation}
Therefore, \eqref{FG} and \eqref{FL} together lead to
\begin{equation}\label{LevyPdf}
    G_{\alpha,1}^{\theta}(x,t) = t^{-1/\alpha}L_{\alpha}^{\theta}(x/t^{1/\alpha}), ~\ -\infty<x<\infty, ~\ t\geq 0.
\end{equation}
Hence, the solution of the space-fractional diffusion equation \eqref{Greenexpr} coincides with the L{\'e}vy strictly stable PDF \eqref{LevyPdf}, in this scenario, \cite{analytic}. For the remainder of this discussion, we restrict ourselves to symmetric distributions, corresponding to $\theta = 0$, and for which $\psi_{\alpha}(k) = |k|^{\alpha}$. 

\subsection{Stable Distributions}
More often than not, one assumes a distribution portraying some sort of real-life phenomenon to be Gaussian. However, opting for this particular approximation is not always the best course of action. In order to set up a given model which behaves more authentically, it is essential that one can incorporate heavier-tailed distributions as well. If one were to pursue a normality assumption, then many real phenomena that can be observed in practice appear near-impossible outcomes of the model, due to how the model was constructed in the first place. The problem lies in the given that the light Gaussian tails make it almost impossible for said events to occur. In other words, the corresponding probability with which they take place is so low that the frequency with which they can (in theory) be observed appears negligibly small. Nevertheless, one still wishes to accurately depict the occurrence of these possible (nevertheless unexpected) jumps in an agent-based random walk model, for instance. Since one such event can only be explained by the occurrence of an outlier in the stochastic model, this justifies the use of heavy-tailed distributions, generally speaking, \cite{MultiStable}. 

The stable distributions represent a flexible parametric family of distributions that augment the class of normal distributions. In particular, any Gaussian distribution can re-appear as a stable distribution itself, for a specific combination of parameters. Even though general stable distributions possess lots of additional properties, therefore making it harder to recognise a connection between them and the multivariate normal distributions, their symmetric sub-family consisting of all elliptical stable distributions very much resembles the normal distributions, both in structure and other related properties. It is sufficient to restrict ourselves to the class of multivariate elliptical stable distributions, which are also more easily simulated than in their most general form, \cite{MultiStable}.

\subsubsection{Drawing from Stable Distributions}
A more general collection of distributions replaces the class of normal distributions in the discrete agent-based model for fractional diffusion. In particular, this is where the stable distributions enter the scene. Any stable distribution is completely determined by its characteristic function (i.e. the Fourier transform of the pdf), as given by \cite{github} and \cite{Nolan}
\begin{equation}\label{stabledistr}
    \varphi(t; \alpha, \beta, \sigma, \mu) = e^{it\mu - \{|\sigma t|^{\alpha}(1-i\beta \text{sgn}(t)\Phi(t))\}}, ~\ \text{ with } ~\
    \Phi(t) = \left \{
    \begin{array}{ll}
         - \frac{2}{\pi}\log|t| & \text{ for $\alpha = 1$}  \\
         \text{tan}{\frac{\pi\alpha}{2}} & \text{ for $\alpha \neq 1$}. 
    \end{array}\right. 
\end{equation}
If we set $\beta =0$, $\sigma =1$ and $\mu=0$ in the characteristic function \eqref{stabledistr} above, then we recover
\begin{equation}\label{FourAlpha}
    \varphi(k; \alpha, \beta = 0, \sigma =1, \mu =0) = e^{-|k|} = \mathcal{F}[L_{\alpha}^0](k),
\end{equation}
which equals the Fourier transform of the L{\'e}vy strictly stable density with parameter $\alpha$ and $\theta = 0$. We call the corresponding L{\'e}vy distribution $L_{\alpha}^0$ a \textit{symmetric $\alpha$-stable distribution}, \cite{github}, \cite{Nolan}, and \cite{Benson2000}. 
%The above information was taken from \cite{github} and is based on \cite{Nolan}.

Since the solution to the space-fractional diffusion equation is of the form \eqref{LevyPdf}, we are actually interested in the scaled random variable $Y=t^{1/\alpha}X$, where $X$ follows a symmetric $\alpha$-stable distribution, as we will denote by $X\sim\text{Stable}(\alpha)$.
In order to draw from the r.v. $Y$ having probability density function \eqref{LevyPdf}, we will first draw a sample from $X$ before scaling the resulting number by the factor $t^{1/\alpha}$, in order to get a randomly drawn sample from $Y$. 

%--------
%\cite{MultiStable}
We are currently only interested in $\alpha$-stable distributions, for which the parameter $\beta$ does not contribute to \eqref{stabledistr}. For general values of $\sigma$ and $\mu$, characteristic functions of the form 
\begin{equation}\label{elliptic}
    \varphi(t; \alpha, \beta=0, \sigma, \mu) = e^{it\mu - (\sigma |t|)^{\alpha}} ~\ \text{ with }\alpha \in (0,2], ~\ \mu \in (-\infty, \infty)  ~\ \text{ and } ~\ \sigma>0
\end{equation}
stem from a symmetric (due to the fact that the skewness parameter $\beta = 0$) subfamily of the stable distributions, collectively known as \textit{elliptical stable distributions}. We dealt with such univariate distributions in the last paragraph and will come to generalize them to multivariate distributions in this paragraph. One recognises that the characteristic function, translated to the multivariate setting, is represented by \cite{MultiStable}
\begin{equation}
    \Psi\left(\Vec{t}~\right) = e^{i \Vec{t}~^T \bm\mu}e^{-|\Vec{t}~^{T}\bm Q\Vec{t}~|^{\alpha/2}}.\label{charmulti}
\end{equation}
\textit{Multivariate elliptic stable distributions} are alternatively called \textit{subgaussian distributions} due to their resemblance to multivariate normal distributions. In fact, a multivariate normal distribution is recovered when setting $\alpha =2$ in \eqref{charmulti}. 
If $\alpha>1$, then $\bm{\mu}$ can still be interpreted as a mean vector. In general, $\bm{Q}$ is a positive definite matrix that determines the dependence structure, it is equal to the covariance matrix of a multivariate normal distribution in case $\alpha=2$, \cite{MultiStable}.

%\paragraph{Simulations}
In order to simulate an agent-based random walk that follows a fractional diffusion model, we will need to draw from a (two-dimensional) multivariate stable distribution. 
Let
\begin{equation}\label{multistable}
    \bm{Z} \sim \text{Stable}\left(\alpha, ~ \bm{Q}, ~ \bm\mu\right)
\end{equation}
denote a multivariate stable random vector having characteristic function coinciding with \eqref{charmulti}.   
For our purposes, where the random vector represents an agent's spatial shift in position coordinates between subsequent time steps, the two-dimensional version of the symmetric $\alpha$-stable distribution $L_{\alpha}^0$---characterised by \eqref{FourAlpha} in the univariate case---is attained by setting $\bm \mu = \bm 0$ and $\bm{Q} = \bm{\text{I}_2}$ in \eqref{multistable}, \cite{MultiStable}.

Let us therefore assume 
\begin{equation}\label{Zstabledim2}
    \bm{Z} \sim \text{Stable}\left(\alpha, ~\ \bm{Q} = \bm{\text{I}_2}, ~\ \bm\mu = \bm{0}\right)
\end{equation}
and suppose we wish to draw a random vector $\bm z$ from the distribution of $\bm Z$. 
As it turns out, simulating a sample from $\bm Z$ is equivalent to simulating a sample $s$ say from the univariate stable distribution
\begin{equation}\label{Sstabledim1}
    S \sim \text{Stable}\left(\frac{\alpha}{2}, ~\ \beta = 1, ~\ \sigma = \left(\cos(\frac{\pi \alpha}{4})\right)^{\frac{2}{\alpha}}, ~\ \mu = 0\right)
\end{equation}
as well as drawing the vector $\bm g$ from the multivariate normal distribution
\begin{equation}\label{Gnormdim2}
    \bm{G} \sim \mathcal{N}( \bm{0}, ~\ \bm{Q} ).
\end{equation}
As pointed out by \cite{MultiStable}, there exists an equivalence of random variables
\begin{equation*}
    \bm Z = \sqrt{S}\bm{G}, ~\ \text{ so that setting } ~\
    \bm z = \sqrt{s}\bm{g} 
\end{equation*}
provides us with a sample of the multivariate stable distribution \eqref{multistable} in the end, \cite{MultiStable}. 

%==============================================================
\section{Numerical Simulations}\label{sec:numerics}
Let us now put theory to the test. In this section, we will simulate simple scenarios of space-fractional diffusion, both on an individualistic level and in a macroscopic setting. We will propose a numerical scheme to discretise the macroscopic diffusion equation and implement it in \texttt{Julia}, from the ground up. Alternative approaches have been undertaken by \cite{Baeumer2008} or \cite{Carrillo}, for instance.
All implementations and simulations take place in \texttt{Julia}, \cite{Juliabezanson2017}. \texttt{Julia} offers the \texttt{StableDistributions.jl} package as part of its coding language. We will use it to draw from stable distributions of the form \eqref{stabledistr}, in particular symmetric $\alpha$-stable distributions as characterised by \eqref{FourAlpha}.
%General stable distributions having characteristic function of the form \eqref{stabledistr} are available under the command ``$\text{Stable}(\alpha, \beta, \sigma, \mu)$''. The command ``$\text{Stable}(\alpha)$'' is shorthand for ``$\text{Stable}(\alpha, \beta=0, \sigma=1, \mu=0)$'' and coincides with a

\subsection{Agent-Based Fractional Diffusion}\label{sec:agentbased}
Our agent-based fractional diffusion simulations run in two dimensions. Denoting by $\bm{X_{t}^i}$ the two-dimensional position vector of the $i$-th particle at time $t$, the governing equation is given by
\begin{equation}
    \bm{X_{t+\Delta t}^i} = \bm{X_{t}^i}  + \Delta t^{1/\alpha}\cdot \sqrt{S}\bm{G}.
\end{equation}
This is consistent with the distribution \eqref{distrXti}, modelling the change in position coordinates for particle $i$, after a time $\Delta t$ (equal to a symmetric L\'evy strictly stable distribution,  as argued by \eqref{LevyPdf}). Let us also recall that the (two-dimensional) multivariate stable distribution \eqref{Zstabledim2} is equivalent to $\sqrt{S}\bm{G}$, with $S$ the one-dimensional stable distribution \eqref{Sstabledim1} and $\bm{G}$ the multivariate normal distribution \eqref{Gnormdim2}. 

\subsubsection{Simulations}
\paragraph{Initial Condition}
In generating this simulation, we have placed $100$ particle agents %\nkscomm{this is the first instance where "virus" or "virotherapy" appear in the text. we should either replace it with something different, e.g. "particle agents", or introduce virotherapy in the introduction. I would suggest that we rather changed the appearance of "virus"} 
(each for the $\alpha=1.5$ and the $\alpha=2$ simulation) at the origin $(0,0)$ at time $t=0$. 

\paragraph{Parameter Values} 

virion count = 100 \big| $\alpha_{\text{frac}} = 1.5$ vs $\alpha_{\text{reg}} = 2$ \big| $\Delta t = 10^{-8}$  \big| $T=10^{-4}$ \big| $\bm{D}_{\text{reg}}=
 \begin{pmatrix}
     0.02 & 0 \\
     0 & 0.02
 \end{pmatrix}$ \big| $\bm{D}_{\text{frac}}=
 \begin{pmatrix}
     1 & 0 \\
     0 & 1
 \end{pmatrix}$
 
\begin{figure}
 \begin{subfigure}{0.49\textwidth}
     \includegraphics[width=\textwidth, trim={90 0 60 43}, clip]{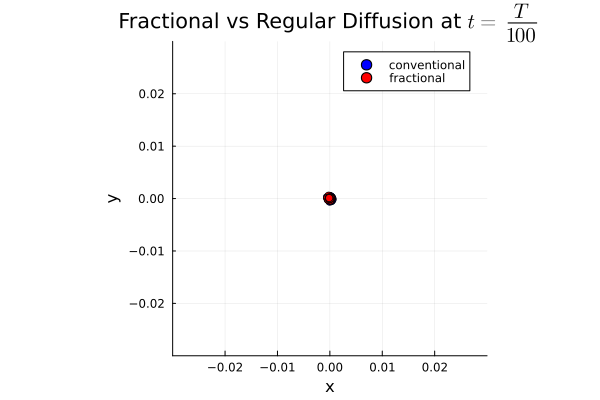}
     \caption{$t = 10^{-6}$}
     \label{micro:a}
 \end{subfigure}
 \hfill
 \begin{subfigure}{0.49\textwidth}
     \includegraphics[width=\textwidth, trim={90 0 60 43}, clip]{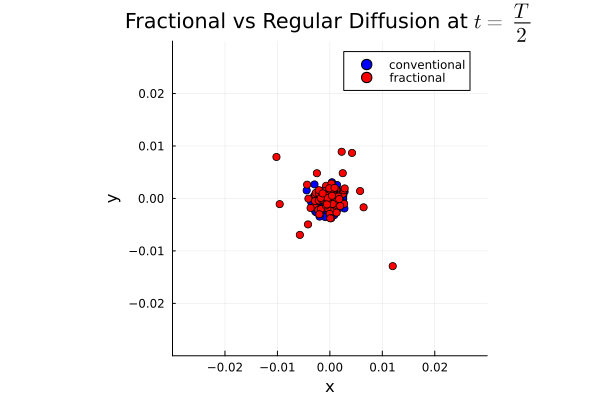}
     \caption{$t = 5\cdot 10^{-5}$}
     \label{micro:b}
 \end{subfigure}
 
 \medskip
 \begin{subfigure}{0.49\textwidth}
     \includegraphics[width=\textwidth, trim={90 0 60 43}, clip]{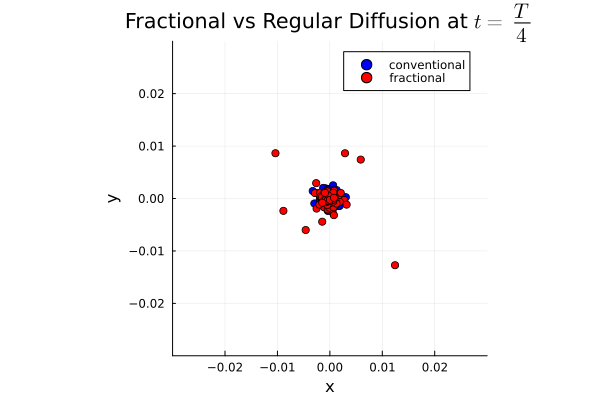}
     \caption{$t = 2.5\cdot 10^{-5}$}
     \label{micro:c}
 \end{subfigure}
 \hfill
 \begin{subfigure}{0.49\textwidth}
     \includegraphics[width=\textwidth, trim={90 0 60 40}, clip]{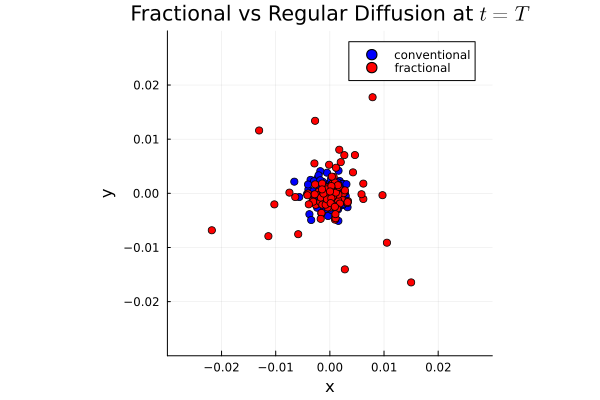}
     \caption{$t = 10^{-4}$}
     \label{micro:d}
 \end{subfigure}
 \caption{Comparison between agent-based (genuinely) fractional and (nearly) classical diffusion. The red dots represent agents dispersing according to the fractional order of $\alpha=1.5$ while the blue agents follow the (nearly) regular ($\alpha =1.99$) diffusion model. Snapshots of the simulation indicate how expansion of the colony evolves at different intermediate time points (a) $t = 10^{-6}$, (b) $t = 5\cdot 10^{-5}$, (c) $t = 2.5\cdot 10^{-5}$, and (d) $t = 10^{-4}$. The respective diffusion coefficients are chosen such that the bulk of the agents is spread over the same convex set in space, in both the fractional and (nearly) regular diffusion case. Nevertheless, the fractional diffusion model exhibits a heavier-tailed distribution that leads to a higher probability of long-range jumps.}
 \label{fig:agentbased}
\end{figure}

\begin{figure}
 \begin{subfigure}{0.49\textwidth}
     \includegraphics[width=\textwidth, trim={25 0 0 20}, clip]{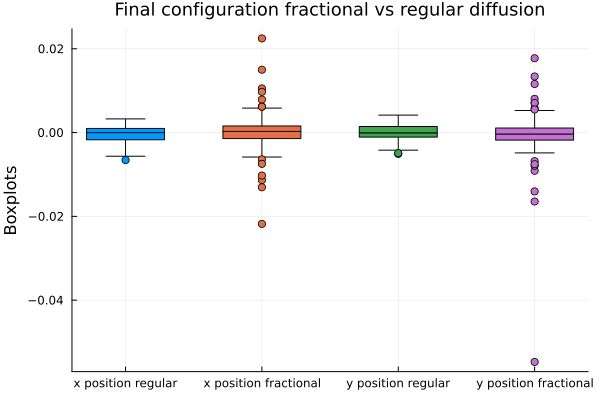}
     \caption{With outliers}
     \label{micro:with}
 \end{subfigure}
 \hfill
 \begin{subfigure}{0.49\textwidth}
     \includegraphics[width=\textwidth, trim={25 0 0 20}, clip]{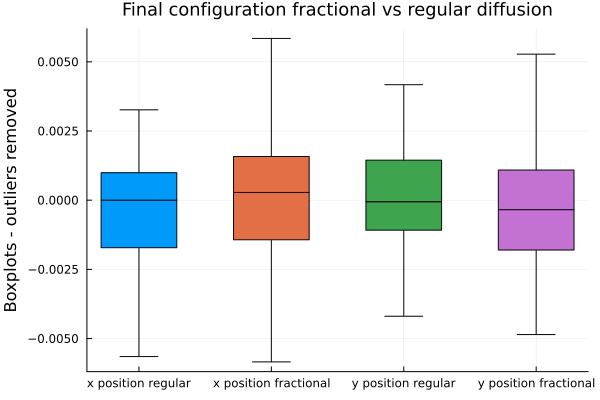}
     \caption{Without outliers}
     \label{micro:without}
 \end{subfigure}
 \caption{Statistical comparison of the agent-based distributions for the regular vs fractional diffusion of Figure \ref{micro:d}, at the final time $T$. (a) Boxplots showing the spread of position coordinates for all agents at the final time $T$, taking the outliers into consideration. (b) Same boxplots after removal of the outliers, for more emphasis on the central distribution. The blue and orange boxplots represent the distribution of the $x$-position coordinates of regularly ($\alpha = 1.99$), respectively fractionally ($\alpha = 1.5$) diffusing agents. Green and purple boxplots indicate how the final $y$-position coordinates are distributed for regularly, respectively (genuinely) fractionally diffusing agents.}
 %\nkscomm{I have trimmed the figures. let's discuss}
 \label{boxplots}
\end{figure}

\paragraph{Results}
Figure \ref{fig:agentbased} exhibits the individual propagation of virion particles initially positioned at the origin. Particles evolving according to the fractional diffusion stochastic diffusion equation for $\alpha =1.5$ are depicted in red. The reason for the different choice of diffusion matrices $\bm{D}_{\text{reg}}$ and $\bm{D}_{\text{frac}}$ lies in us wanting to compare the $\alpha=1.5$ against the classical $\alpha=2$ diffusive behaviour in a systematic way. Since every random jump takes place according to the two ($x$- and $y$-) coordinates, we have opted for a two-dimensional stable versus normal distribution. The $\alpha=1.5$ distribution showcases a heavy tail that can be connected to the rarity of the jumping event. Moreover, since we assume an independence of the coordinates, meaning that the $x$- and $y$- jump lengths are uncorrelated, we opt for diagonal diffusion matrices. The bulk of particles sticking together around $(0,0)$ even after $t=0$ is what we will refer to as the ``clump''. The spread of both blue and red clumps depends on the diffusion matrix entry values. We wanted to showcase an identical diffusion behaviour between the clumps themselves and in order to attain that image, we had to choose different diffusion matrices. Choosing for the fractional $\alpha=1.5$ diffusion model the identity matrix, we calibrated the diffusion matrix $\bm{D}_{\text{reg}}=\begin{pmatrix}
    0.2 & 0 \\ 0 & 0.2
\end{pmatrix}$ for the standard $\alpha =2$ diffusion, through visual inspection of the simulation results. We did so by generating the convex hull of the regularly diffusing particles at their end position $t=T$. Based on this convex set, we visually determined the fractional $\alpha =1.5$ diffusion coefficients needed in order for the red and blue clumps to overlap. As it turns out, this ensured the particle bulks to be similar at all (previous) times, with the exception of some red outliers floating around further away the main clumps. 

In Figure \ref{boxplots}, we have created boxplots of the diffusing particles' end positions at time $t=T$. This provides us with a statistical measure of comparison. Figure \ref{micro:without} shows the inter quartile ranges (with median at $0$) to be near identical upon comparison of the regular and fractional coordinates. The main difference between the fractional and regular boxplots for $x$- and $y$- coordinates lies in the amount of outliers. Excluding outliers, the remaining particles function as clumps of the same diameter having similar qualitative behaviour. Choosing the same diffusion matrices (say the identity matrix for instance) for $\alpha=1.5$ and $\alpha=2$ would cause a difference in the intrinsic speeds of the individual $\alpha=1.5$ versus $\alpha=2$ agents, causing the colonies to move at vastly different rates themselves. Excluding the outliers in the scenario depicted by Figure \ref{fig:agentbased}, what remains of the colony can be said to move at the same speed, however. 
% {\color{red}
% systematic comparison speed, diffusion coefficient, alpha?}
% {\color{blue}Insert simulations same diffusion matrix? Modify figures.}
%(biggest convex set containing all elements)
%other application: small mutations from time to time $\rightarrow$ jump/tails over generations

\subsection{Macroscopic Fractional Diffusion}
For most fractional differential equations, an exact analytical solution cannot be easily obtained. Even though several different approaches, such as integral transformation methods, have become available and can be used to obtain precise solutions in some, special cases, more often than not they will fail to help us. Hence, one has to appeal to numerical methods. Focusing on the case $\beta =1$, we are left with a Riesz space-fractional diffusion equation for which the goal is to find a (higher-order) numerical scheme, based on a suitable discretisation of the Riesz derivative. Multiple individuals have already ventured coming up with numerical schemes to estimate the fractional differential equation's actual solution. However, higher-order numerical schemes that implement mass-conserving boundary conditions seem to be scarce in the literature, \cite{scheme}. The fractional Laplacian, and its boundary conditions in particular have been studied in detail by \cite{Caffarelli2007}. \cite{Cusimano2018} propose a different discretisation of the spectral fractional Laplacian on bounded domains. Their approach is based on the integral formulation of the operator via the heat-semigroup formalism.
%In this work, we propose novel discretizations of the spectral fractional Laplacian on bounded domains based on the integral formulation of the operator via the heat-semigroup formalism.  can be implemented on possibly irregular bounded domains, and can naturally handle different types of boundary constraints. 
Among the available numerical approaches, spectral methods have gained attention for their ability to provide high-accuracy solutions while reducing computational memory requirements, particularly in time-fractional problems (see \cite{Carrillo2023}).
In this paper, we focus on an explicit scheme for solving the fractional superdiffusion equation with corresponding parameters $\beta =1, 1<\alpha<2$. Instead of directly enforcing a fractional equivalent of zero-flux boundary conditions, we are using a slightly different approach to deal with the boundary. Namely, we will implement periodic boundaries in such a way that the scheme becomes symmetric at every spatial grid point. This will be made precise as we develop the numerical scheme. 

\subsubsection{Setup}
We propose a scheme---that is second-order accurate in space and explicit in time---for finding the approximate solution of a space-fractional diffusion equation. Augmented with periodic boundaries, we introduce a novel scheme to approach fractional diffusion. A Taylor expansion method, including the shifted Gr{\"u}nwald-Letnikov operator, lies at the heart of the spatial Riesz derivative's approximation. The Gr{\"u}nwald-Letnikov operator, and other numerical methods for computing fractional derivatives, have been studied in detail by \cite{Podlubny1999}. Let $u(x,t)$ describe the considered particle's density at the point $x$ at time $t$. Note that in a first instance, we will reduce ourselves to a one-dimensional, bounded domain. The defining system of equations is then given by  
\begin{align}
    \pdv{}{t}u(x,t) &= D \pdv[\alpha]{}{|x|}u(x,t), ~\ a<x<b, ~\ 0<t\leq T,  ~\ \text{with}\label{fracdiff}\\
    u(x,0) &= \phi_0(x),\label{fracdiff0}
\end{align}
on a macroscopic level. We further restrict ourselves to a finite time interval $[0,T]$ on which we compute the numerical solution. In the fractional diffusion equation \eqref{fracdiff}, $D$ represents the diffusion coefficient. $\phi_0$ describes a sufficiently smooth initial function, \cite{XieFang}. We will employ the definition \eqref{defRiesz}, characterising the Riesz derivative in terms of left- and right Riemann-Liouville derivatives. 

\subsubsection{Gr{\"u}nwald-Letnikov}
The---weighted and shifted---Gr{\"u}nwald-Letnikov operator will be used to discretise the fractional Riesz derivative $\pdv[|\alpha|]{}{x}$. 
\begin{definition}[Gr{\"u}nwald-Letnikov operator]\label{defgrun}
    The weighted and shifted Gr{\"u}nwald-Letnikov operators are defined by 
    \begin{align}
        { }_{\text{L}}D_{h,p,q}^{\alpha}u(x) &= \frac{\lambda_1}{h^{\alpha}} \sum_{j=0}^{\infty}g_j^{(\alpha)}u(x-(j-p)h) + \frac{\lambda_2}{h^{\alpha}}\sum_{j=0}^{\infty}g_j^{(\alpha)}u(x-(j-q)h) ~\ \text{ and }\\
        { }_{\text{R}}D_{h,p,q}^{\alpha}u(x) &= \frac{\lambda_1}{h^{\alpha}} \sum_{j=0}^{\infty}g_j^{(\alpha)}u(x+(j-p)h) + \frac{\lambda_2}{h^{\alpha}}\sum_{j=0}^{\infty}g_j^{(\alpha)}u(x+(j-q)h),
    \end{align}
    for integers $p,q$ with $p\neq q$, where $\lambda_1= \frac{\alpha-2q}{2(p-q)}$ and with $\lambda_2= \frac{2p-\alpha}{2(p-q)}$. The fractional binomial coefficients appear as 
    \begin{equation}\label{coeffgrun}
        g_j^{(\alpha)} = (-1)^{j}\frac{\Gamma(\alpha+1)}{\Gamma(\alpha-j+1)\Gamma(j+1)}. 
    \end{equation}
\end{definition}
The next result explains how these operators can be used to yield second-order accurate approximations in space, for left- and right Riemann-Liouville derivatives respectively. The following lemma is taken from (Lemma 1, \cite{XieFang}) and formulated here in a slightly altered form.
\begin{lemma}\label{lemgrun}
    Let $u(x)\in L^1(\R )$ and suppose that ${ }_{\text{LRL}}D_{a,x}^{\alpha+2}u(x)$ and ${ }_{\text{RRL}} D_{x,b}^{\alpha+2} u(x)$---along with their Fourier transforms---belong to $L^1(\R )$. Then the following second-order approximations are obtained, namely
    \begin{align*}
       { }_{\text{L}}D_{h,p,q}^{\alpha}u(x)= { }_{\text{LRL}}D_{a,x}^{\alpha}u(x) + \mathcal{O}(h^2) ~\ \text{ and } ~\
        { }_{\text{R}}D_{h,p,q}^{\alpha}u(x)={ }_{\text{RRL}} D_{x,b}^{\alpha} u(x) + \mathcal{O}(h^2),
    \end{align*}
    and this holds uniformly for all $x\in \R $. 
\end{lemma}
The next subsection showcases how this lemma can be deployed to develop a numerical scheme for the fractional Riesz-diffusion equation. 

\subsubsection{Numerical Scheme for the Macroscopic Equation}
Our numerical grid consists of a discrete series of space- and time tuples $(x_i,t_n)$, where 
\begin{equation}
    x_i = a+ih, ~\ i=0, \hdots, M ~\ \text{ and } ~\ t_n = n\tau, ~\ n = 0,\hdots, N,
\end{equation}
for certain $N,M \in \mathbb{N}_{>0}$. 
The uniform space- and time step sizes $h$ and $\tau$ are defined by
\begin{equation}
    h=\frac{(b-a)}{M} ~\ \text{ and } \tau = \frac{T}{N}. 
\end{equation}
For simplicity, note that we here adopt a uniform grid.
We further denote
\begin{equation*}
    u_i^n = u(x_i,t_n) ~\ \text{ for short, with $i=0, \hdots, M$ and $n=1, \hdots, N$.}
\end{equation*}

In order to come up with a numerical scheme solving \eqref{fracdiff}, we will first discretise the Riesz derivative using the weighted and shifted Gr{\"u}nwald-Letnikov operator. As an immediate consequence of the previous Lemma \ref{lemgrun}, together with Definition \ref{defgrun} of the weighted and shifted Gr{\"u}nwald-Letnikov operators, we can use the second-order approximations 
\begin{align*}
     { }_{\text{LRL}}D_{a,x}^{\alpha}u(x,t) &=\frac{\lambda_1}{h^{\alpha}} \sum_{j=0}^{\infty}g_j{(\alpha)}u(x-(j-p)h,t) + \frac{\lambda_2}{h^{\alpha}}\sum_{j=0}^{\infty}g_j^{(\alpha)}u(x-(j-q)h,t) + \mathcal{O}(h^2) ~\     \text{ and }
    \\
    { }_{\text{RRL}} D_{x,b}^{\alpha} u(x,t) &= \frac{\lambda_1}{h^{\alpha}} \sum_{j=0}^{\infty}g_j{(\alpha)}u(x+(j-p)h,t) + \frac{\lambda_2}{h^{\alpha}}\sum_{j=0}^{\infty}g_j^{(\alpha)}u(x+(j-q)h,t) + \mathcal{O}(h^2)
\end{align*}
to do so. 
Note that there are many other approximations of the Riesz derivative at hand, one of them being the \textit{fractional central difference method}. See \cite{scheme} for more information on this topic (as well as an even higher-order $\mathcal{O}(\tau^2 + h^4)$ numerical scheme proposed by the authors Ding et al. themselves).

Since Definition \ref{defgrun} allows us to choose the values of $p$ and $q$ freely, we set $p=1$ and $q =0$ to obtain the simplified values of $\lambda_1 = \frac{\alpha}{2}$ and $\lambda_2= \frac{2-\alpha}{2}$ for the coefficients. Moreover, seeing as we only consider a bounded domain (on which $u$ is defined, i.e. non-zero in particular), the infinite Gr{\"u}nwald-Letnikov sums reduce to finite sums. We find
\begin{align*}
     \left.{ }_{\text{LRL}}D_{a,x}^{\alpha}u(x,t)\right\vert_{(x_i,t_n)} &=\frac{\lambda_1}{h^{\alpha}} \sum_{j=0}^{\infty}g_j^{(\alpha)}u^n_{i-j+1} + \frac{\lambda_2}{h^{\alpha}}\sum_{j=0}^{\infty}g_j^{(\alpha)}u^n_{i-j} + \mathcal{O}(h^2)\\
    &=\frac{\lambda_1}{h^{\alpha}} \sum_{j=0}^{i+1}g_j^{(\alpha)}u^n_{i-j+1} + \frac{\lambda_2}{h^{\alpha}}\sum_{j=0}^{i}g_j^{(\alpha)}u^n_{i-j} + \mathcal{O}(h^2)\\
    &=\frac{\lambda_1}{h^{\alpha}} \sum_{j=0}^{i+1}g_j^{(\alpha)}u^n_{i-j+1} + \frac{\lambda_2}{h^{\alpha}}\sum_{j=1}^{i+1}g_{j-1}^{(\alpha)}u^n_{i-j+1} + \mathcal{O}(h^2),
\end{align*}
due to the fact that $u^n_k$ is only defined for integers $k$ for which $0\leq k \leq M$. The latter criterion then restricts the upper bound for the summation index $j$ (in terms of $i$) in the infinite Gr{\"u}nwald-Letnikov sums, and renders them finite.
We ultimately arrive at
\begin{align}
     \left.{ }_{\text{LRL}}D_{a,x}^{\alpha}u(x,t)\right\vert_{(x_i,t_n)} 
    &=\frac{\lambda_1}{h^{\alpha}} \sum_{j=0}^{i+1}g_j^{(\alpha)}u^n_{i-j+1} + \frac{\lambda_2}{h^{\alpha}}\sum_{j=1}^{i+1}g_{j-1}^{(\alpha)}u^n_{i-j+1} + \mathcal{O}(h^2)\nonumber\\
    &=\frac{1}{h^{\alpha}} \sum_{j=0}^{i+1} w_j^{(\alpha)}u^n_{i-j+1}  + \mathcal{O}(h^2), \label{LRLGrun}
\end{align}
where (with the Gr{\"u}nwald weights defined as in \eqref{coeffgrun}, recalling the explicit values for $\lambda_1$, $\lambda_2$)
\begin{align}
    w_j^{(\alpha)}&= \left \{ 
    \begin{array}{ll}
         \frac{\alpha}{2} g_0^{(\alpha)}& \text{if $j=0$,}   \\
        \frac{\alpha}{2}g_j^{(\alpha)}+\frac{2-\alpha}{2}g_{j-1}^{(\alpha)} & \text{elsewhere, for $j=1, \hdots, M$.}
    \end{array}
    \right. 
\end{align}
Similarly, using the Gr{\"u}nwald-Letnikov operator to approximate the right Riemann-Liouville derivative on a bounded domain, our calculations lead to
\begin{align}
    \left.{ }_{\text{RRL}} D_{x,b}^{\alpha} u(x,t)\right\vert_{(x_i,t_n)} &= \frac{\lambda_1}{h^{\alpha}} \sum_{j=0}^{\infty}g_j^{(\alpha)}u^n_{i+j-1} + \frac{\lambda_2}{h^{\alpha}}\sum_{j=0}^{\infty}g_j^{(\alpha)}u^n_{i+j} + \mathcal{O}(h^2)\nonumber\\
    &=\frac{\lambda_1}{h^{\alpha}} \sum_{j=0}^{M-i+1}g_j^{(\alpha)}u^n_{i+j-1} + \frac{\lambda_2}{h^{\alpha}}\sum_{j=0}^{M-i}g_j^{(\alpha)}u^n_{i+j} + \mathcal{O}(h^2)\nonumber\\
    &=\frac{\lambda_1}{h^{\alpha}} \sum_{j=0}^{M-i+1}g_j^{(\alpha)}u^n_{i+j-1} + \frac{\lambda_2}{h^{\alpha}}\sum_{j=1}^{M-i+1}g_{j-1}^{(\alpha)}u^n_{i+j-1} + \mathcal{O}(h^2)\nonumber\\
    &=\frac{1}{h^{\alpha}}\sum_{j=0}^{M-i+1}w_j^{(\alpha)}u^n_{i+j-1} + \mathcal{O}(h^2). \label{RRLGrun}
\end{align}
Recalling the definition \eqref{defRiesz} of the Riesz derivative in terms of (left- and right) Riemann-Liouville derivatives 
\begin{align*}
    \pdv[\alpha]{}{|x|}u(x,t) &= -c_{\alpha} \left( { }_{\text{LRL}}D_{a,x}^{\alpha} + { }_{\text{RRL}} D_{x,b}^{\alpha}\right) u(x,t),
\end{align*}
\eqref{LRLGrun} and \eqref{RRLGrun} together lead to the second order approximation
\begin{align}\label{secondorder}
    \left.\pdv[\alpha]{}{|x|}u(x,t)\right\vert_{(x_i,t_n)}
    &= -\frac{c_{\alpha}}{h^{\alpha}} \left(\sum_{j=0}^{i+1} w_j^{(\alpha)}u^n_{i-j+1} + \sum_{j=0}^{M-i+1} w_j^{(\alpha)}u^n_{i+j-1}\right) + \mathcal{O}(h^2),
\end{align}
for $i = 1, \hdots, M-1$.

%\paragraph{Explicit periodic scheme \cite{scheme}}
We still need to develop a numerical scheme which computes $u$ as defined by \eqref{fracdiff} at every grid point (not just at the internal points which discard the boundary of the spatial grid). 
First of all, we replace the Riesz fractional derivative by its second order approximation \eqref{secondorder}.
As for the ordinary time derivative, Taylor expansion of the function $u$ leads to
\begin{align}\label{tayloru}
    u_i^{n-1} = u_i^n - \tau \pdv{u_i^n}{t} + \mathcal{O}(\tau^2).
\end{align}
and brings us to the first order forward Euler difference
\begin{align}\label{tayloruexpl}
    \pdv{u_i^n}{t} \approx \frac{u_i^n - u_i^{n-1}}{\tau}.
\end{align}

Combining \eqref{secondorder} and \eqref{tayloruexpl}, we come to the $\mathcal{O}(\tau, h^2)$-approximation 
\begin{align}\label{numericexpl}
u_i^n &= u_i^{n-1} - \tau D \frac{c_{\alpha}}{h^{\alpha}}\left \{\sum_{j=0}^{i+1} w_j^{(\alpha)}u^{n-1}_{i-j+1} + \sum_{j=0}^{M-i+1} w_j^{(\alpha)}u^{n-1}_{i+j-1}\right\} + \mathcal{O}(\tau, h^2)
\end{align}
for $i=1,\hdots, M-1$ and $n=0,\hdots,N$, in the first instance. 
Let $U_i^n$ be the numerical solution analogue of $u_i^n$. Then, getting rid of the $\mathcal{O}(\tau, h^2)$-term in \eqref{numericexpl}, we arrive at the explicit-in-time scheme
\begin{align}\label{schemeeq}
    U_i^n 
    &= 
    U_i^{n-1} - r \sum_{j=0}^{i+1} w_j^{(\alpha)}U_{i-j+1}^{n-1} - r \sum_{j=0}^{M-i+1}w_j^{(\alpha)}U_{i+j-1}^{n-1},
\end{align}
where we have substituted the coefficient
\begin{align}\label{r}
    r&=\frac{D\tau c_{\alpha}}{h^{\alpha}}.
\end{align}
Note that this scheme is only valid for indices $1 \leq i \leq M-1$ (due to the approximation of the spatial derivative resulting from the Gr{\"u}nwald-Letnikov operator). 

In order to arrive at a more compact, yet still comprehensive, configuration of the scheme, the goal is to rewrite the discrete equations in the form of a vector inner product.
%The first step in doing so is by changing the index $j$ within the summations of \eqref{schemeeq}.
% In the first sum on either side of the equality sign, let us put
% $i-j+1 = k \text{ i.e. } i-k+1 = j$.
% Since
% $j = 0, \hdots, i+1,  \text{ we find }
%     -j= -i-1, \hdots, 0  \text{ hence }
%     k = i-j+1 = 0, \hdots, i+1$.
% Next, let us also change the index of the second sum. Here, we set
% $i+j-1 = l \text{ i.e. }
%     l-i+1 = j$.
% Since
% $j = 0, \hdots, M-i+1,  \text{ we find }
%     l = i+j-1 = i-1, \hdots M$.
Applying two index changes in \eqref{schemeeq},
\begin{align*}
    U_i^n
    &= 
    U_i^{n-1} - r \sum_{k=0}^{i+1} w_{i-k+1}^{(\alpha)}U_{k}^{n-1} - r \sum_{l=i-1}^{M}w_{l-i+1}^{(\alpha)}U_{l}^{n-1}
\end{align*}
is what we are left with, for $i=1,\hdots, M-1$ and $n=1,\hdots,N$. 

Notice how the numerical scheme \eqref{numericexpl}
\begin{align}\label{explschemeM2}
        U_{M_2}^{n} &= U_{M_2}^{n-1} - \tau D \frac{c_{\alpha}}{h^{\alpha}}\left \{\sum_{k=0}^{M_2+1}w_{M_2-k+1}^{(\alpha)}U_k^{n-1}+\sum_{l=M_2-1}^M w_{l-M_2+1}^{(\alpha)}U_l^{n-1}\right\}  +\mathcal{O}(\tau, h^2)
\end{align}
encapsulates symmetry at the middle point $i=\frac{M}{2}=:M_2$. Indeed, we can rewrite \eqref{explschemeM2} in terms of an inner product, namely
\begin{align*}
    U_{M_2}^n = \left\langle \overrightarrow{w}, \overrightarrow{U^{n-1}} \right\rangle ~\ \text{ with } ~\ \overrightarrow{w}
    &=
    \begin{pmatrix}
    -rw^{(\alpha)}_{M_2+1},
    -rw^{(\alpha)}_{M_2},
    \hdots,
    -rw_4^{(\alpha)},
    -rw_3^{(\alpha)},
    -rw_0^{(\alpha)}-rw_2^{(\alpha)},\\
    1-2rw_0^{(\alpha)},\\
    -rw_0^{(\alpha)}-rw_2^{(\alpha)},
    -rw_3^{(\alpha)},
    -rw_4^{(\alpha)},
    \hdots,
    -rw^{(\alpha)}_{M_2},
    -rw^{(\alpha)}_{M_2+1}
\end{pmatrix}^T
\end{align*}
%
% \begin{align}
%     \overrightarrow{w}&=
%     \begin{pmatrix}
%     -rw^{(\alpha)}_{M_2+1}, 
%     -rw^{(\alpha)}_{M_2},
%     \cdots,\\
%     -rw_4^{(\alpha)},
%     -rw_3^{(\alpha)},
%     -rw_0^{(\alpha)}-rw_2^{(\alpha)},
%     1-2rw_0^{(\alpha)},
%     -rw_0^{(\alpha)}-rw_2^{(\alpha)},
%     -rw_3^{(\alpha)},
%     -rw_4^{(\alpha)},\\
%     \cdots,
%     -rw^{(\alpha)}_{M_2},
%     -rw^{(\alpha)}_{M_2+1}
% \end{pmatrix}^T
% \end{align}
% 
a vector of length $M+1$ that carries the coefficients of the sums in \eqref{explschemeM2}. The vector $\overrightarrow{w}$ can be seen to act as some kind of weight vector, according to which the different entries of $\overrightarrow{U^{n-1}}$ should be taken into account in the sum. Note that the vector $\overrightarrow{w}$ itself is symmetric around its middle $(M_2+1)$-th row (with a corresponding value of $1-2rw_0^{(\alpha)}$).

We would like to retain this symmetry for the remaining grid points $i\neq M_2$. Essentially, we want to use equal amounts of information to the left and right from the point of consideration. We will do so by folding the domain closed (think about the spatial domain as being a circle rather than a straight line). The point at index $i$ will act as the new middle point and we will make up for the asymmetry by splitting the domain in half with respect to this point $x_i = a+ih$ (whilst still visualising the domain as being circular). Here, $x_M$ and $x_0$ are considered to be neighbouring points closing the circle.

The weight vector $\overrightarrow{w}$ remains a constant vector. The only thing we need to change in the numerical scheme at every index $i$ is the order of the grid points lining up with the corresponding---constant---weights. We do so by introducing a permutation matrix $P_{i,M}$ that depends on both the index of consideration $i$ as well as the index $M$ of the last grid point.
In other words, for $i\neq M_2$, one writes 
\begin{align}\label{numschememac}
    U_{i}^n = \left\langle \overrightarrow{w}, P_{i,M}\overrightarrow{U^{n-1}} \right\rangle
~\ \text{ where } ~\
    \overrightarrow{U^{n-1}} = (U_i^{n-1})_{i=0,\hdots,M} = (U_0^{n-1}, \hdots, U_M^{n-1})^T.
\end{align}
The inner product $\langle \cdot , \cdot \rangle$ represents the standard dot product.

For $i=0, \hdots, M$, let us introduce the standard basis vectors $e_{i,M}$ of $\R ^{M+1}$ by setting
\begin{align*}
    e_{i,M} = (\delta_{i,j})_{j=0,\hdots,M}.
\end{align*}
Using the above notation, we can recover the $(M+1)\times(M+1)$ identity matrix 
\begin{align*}
    \text{Id}_{(M+1)\times(M+1)} = \begin{pmatrix}
        e_{0,M}^T, e_{1,M}^T, \hdots, e_{M,M}^T
    \end{pmatrix}^T.
\end{align*}
This makes it easy to define the permutation matrix $P_{i,M}$. For $0\leq i < M_2$, the $k$-th row ($k=1, \hdots, M+1$) of $P_{i,M}$ is defined through
\begin{align*}
    (P_{i,M})[k,:] &=\left \{ \begin{array}{ll}
        e_{M-(M_2-i)+k,M} &  \text{ for $k= 1, \hdots,(M_2-i)$}\\
        e_{k-(M_2-i+1),M} & \text{ for $k= (M_2-i+1),\hdots,(M+1)$}.
    \end{array}\right.
\end{align*}
If $i=M_2$, we recover the identity matrix
\begin{align*}
    P_{M_2,M} &=\text{Id}_{(M+1)\times(M+1)},
\end{align*}
while for $M_2<i\leq M$
\begin{align*}
    (P_{i,M})[k,:] &=\left \{ \begin{array}{ll}
        e_{(i-M_2)+k-1,M} &  \text{ for $k= 1, \hdots,(M-i+M_2+1)$}\\
        e_{k-(M-i+M_2+2),M} & \text{ for $k= (M-i+M_2+2),\hdots,(M+1)$}.
    \end{array}\right.
\end{align*}

\subsubsection{Simulations for the Macroscopic Equation}
We will now use the numerical scheme \eqref{numschememac} to simulate fractional diffusion on a macroscopic scale, in \texttt{Julia}. 
\paragraph{Initial Condition}
We work with a smooth initial condition that is symmetric around $x=0$.
Let us take for $u(x,0)$ the following (unnormalized) Gaussian density with standard deviation $\sigma = 0.2$, centered around $0$
\begin{align*}
    u_{\text{in}}(x,0) &=  e^{-x^2/(2 \sigma^2)}.
\end{align*}
We will further choose $a=-b$ to enforce symmetry of the domain around $0$. 
%continuous/smooth initial condition
%"simulation_results_3_301_200_0_2.jld2"

\paragraph{Chosen Parameter Values}
Unless specified otherwise, the values we are working with are given by
$\sigma=0.2$ \big| $M = 300$ \big| $N=199$  \big|$(a,b)=(-3,3)$ \big| $T=0.005$  \big| $D_{\text{reg}}=0.02$  \big| $D_{\text{frac}}=1$ \big| $\tau = 2.5125628140703518\cdot10^{-5}$ \big| $h = 2 \cdot 10^{-2}$.
\begin{figure}
 \begin{subfigure}{0.49\textwidth}
     \includegraphics[width=\textwidth,trim={25 0 0 23}, clip]{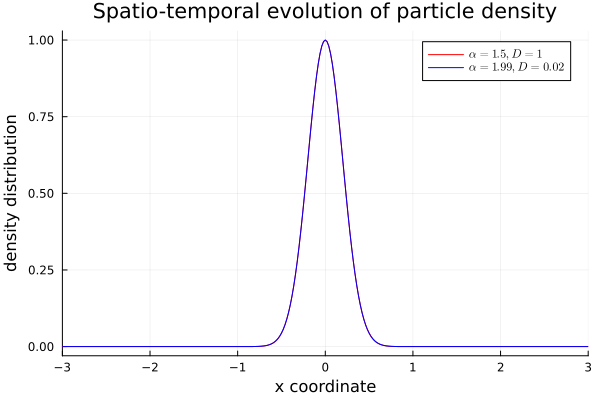}
     \caption{$t = 0$}
     \label{macro:a}
 \end{subfigure}
 \hfill
 \begin{subfigure}{0.49\textwidth}
     \includegraphics[width=\textwidth,trim={25 0 0 23}, clip]{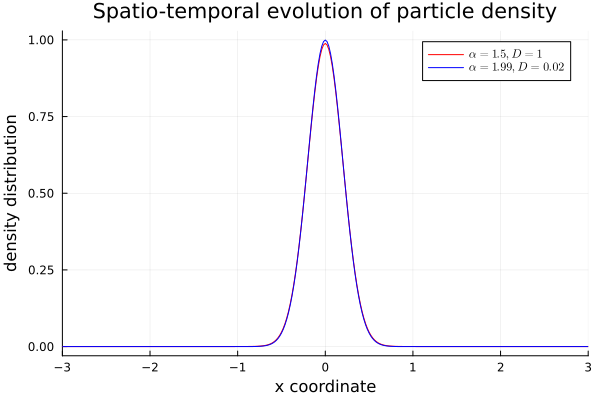}
     \caption{$t = T/4$}
     \label{macro:b}
 \end{subfigure}
 
 \medskip
 \begin{subfigure}{0.49\textwidth}
     \includegraphics[width=\textwidth, trim={25 0 0 23}, clip]{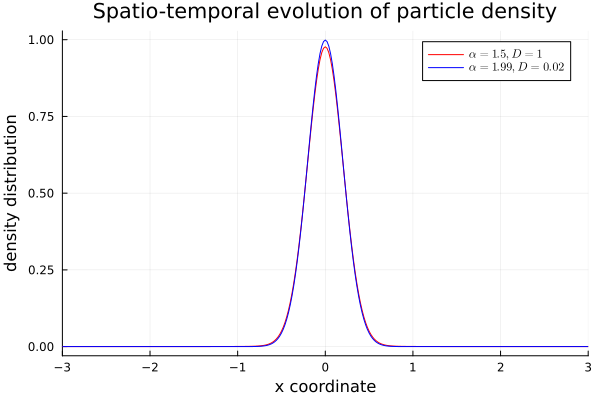}
     \caption{$t = T/2$}
     \label{macro:c}
 \end{subfigure}
 \hfill
 \begin{subfigure}{0.49\textwidth}
     \includegraphics[width=\textwidth, trim={25 0 0 23}, clip]{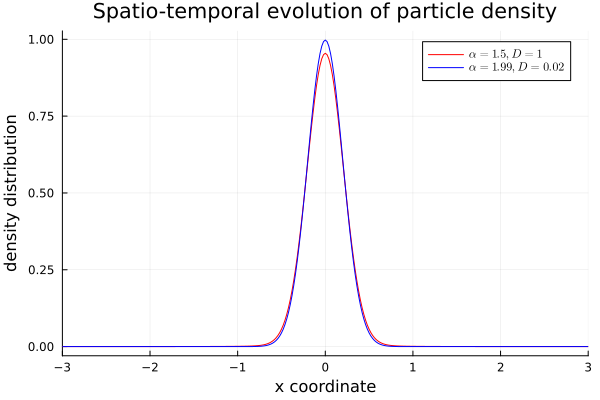}
     \caption{$t = T$}
     \label{macro:d}
 \end{subfigure}

\caption{Macroscopic fractional vs regular diffusion, simulation results in one spatial dimension. The red curve represents the density spread of a fractionally ($\alpha =1.5$) diffusing population with diffusion coefficient $D_{\text{frac}} = 1$. The blue curve ($\alpha = 1.99$) indicates the density evolution of a regularly diffusing density, corresponding to a diffusion coefficient of $D_{\text{reg}}=0.02$. The diffusion coefficients were chosen such that the bulk of the densities remain co-localised over time. A pairwise comparison between regular and genuinely fractional diffusion at different time points (a) $t=0$, (b) $t=T/4$, (c) $t=T/2$, and (d) $t=T$ is presented in the above Figure \ref{fig:macro}.}
\label{fig:macro}
\end{figure}

\begin{figure}
\begin{subfigure}{0.49\textwidth}
    \includegraphics[width=\textwidth, trim={25 0 0 25}, clip]{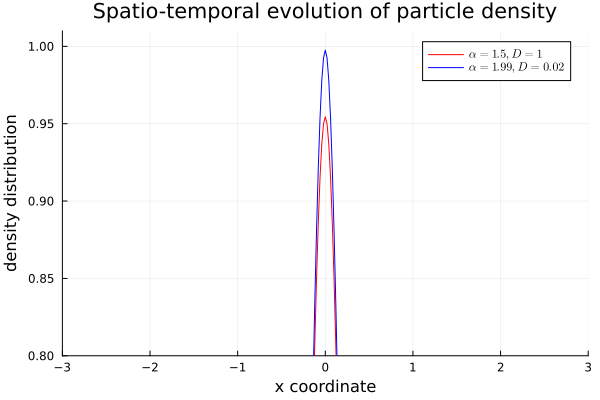}
    \caption{Peak at $t=T$}
    \label{zoom:peak}
\end{subfigure}
\hfill
\begin{subfigure}{0.49\textwidth}
    \includegraphics[width=\textwidth, trim={25 0 0 25}, clip]{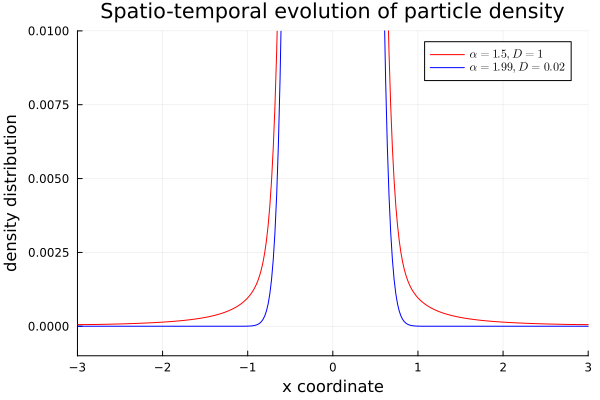}
    \caption{Tails at $t=T$}
    \label{zoom:tails}
\end{subfigure}
\caption{Macroscopic fractional vs regular diffusion simulation: zoomed-in picture at time $t=T$. To display the distinction between regularly ($\alpha=1.99$, blue curve) and genuinely fractional ($\alpha=1.5$, red curve) diffusion, we focus on the density (a) peaks and (b) tails, at the final time of the macroscopic density simulations, by zooming in on Figure \ref{macro:d}.}
\label{zoom_macro}
\end{figure}

% \paragraph{Results}\nkscomm{No need to have paragraph here; we just proceed with the text}
Figure \ref{fig:macro} showcases a pairwise comparison between the $\alpha =1.5$ and $\alpha =1.99$ density curves on the fixed domain $(-3,3)$. Note that we have explicitly chosen $\alpha=1.99$ instead of $2$ to abide by the condition $1<\alpha<2$ that was imposed on the numerical scheme. We have included snapshots of the density evolution, taken at different time points in the interval $[0,T]$. The heavy-tailedness inherent to fractional diffusion is also visible on the macroscopic level. As time evolves, the red curve (corresponding to $\alpha =1.5$) tends to spread out more in space. The (one-dimensional) diffusion coefficients $D_{\text{reg}}=0.02$ and $D_{\text{frac}}=1$ were chosen in accordance with the agent-based diffusion model simulations of Section \ref{sec:agentbased}. 

The numerical scheme used in this study is inherently periodic. While a rigorous analysis of its mass conservation- , positivity-preserving- , consistency- , and convergence properties is necessary to fully assess its accuracy and stability, this study will have to be postponed for future work. Of particular importance are the positivity and mass conservation of our scheme, in order for it to be physically relevant. Our preliminary numerical investigations indicate that the proposed scheme exhibits some mass dissipation, which may be influenced by boundary effects, and certainly by numerical discretisation. However, the use of this numerical scheme is justified by the fact that the macroscopic results, shown in Figure \ref{fig:macro}, are consistent with the atomistic-level observations of Figure \ref{fig:agentbased}. Future studies will have to focus on quantifying and mitigating this mass dissipation, in order to enhance accuracy and reliability of the scheme. 

\section{Conclusion and Outlook}

This paper discusses the role of fractional calculus in the modelling of diffusion processes, at both the micro- and macroscopic level. The derivation of fractional reaction-diffusion equations from CTRWs demonstrates that fractional derivatives arise in a natural way, in systems where waiting times and/or jump lengths follow heavy-tailed distributions. More specifically, both standard and anomalous diffusion equations emerge as a result of taking the long-time, large-scale limit of a CTRW---after applying appropriate rescaling procedures. An infinite variance in jumps leads to the fractional space derivative, while infinite mean waiting times result in fractional Caputo derivatives, \cite{SubdiffFD}. This distinction is especially relevant when modelling transport phenomena in heterogeneous media, where fractional derivatives become essential to accurately describe the dynamics of the system.

At the microscopic level, we have considered fractional stochastic differential equations. By including L\'evy flights, we were able to generalise classical stochastic processes. Instead of assuming Gaussian-distributed steps, we have simulated random walks by drawing from stable distributions, effectively capturing non-local transport phenomena. At the macroscopic level, we have developed a numerical scheme for solving Riesz space-fractional diffusion equations augmented with periodic boundary conditions. Two-sided fractional diffusion equations are becoming increasingly prominent in many applications, with symmetric diffusion on a bounded domain giving rise to the Riesz space-fractional derivative, \cite{BCtwosided}. Driven by the scarcity of numerical methods that preserve mass in case of fractional diffusion, \cite{XieFang}, we have explored a practical alternative for discretising the Riesz fractional diffusion equation.

Our numerical experiments further reinforce the distinctive behaviour between fractional diffusion and its classical counterpart. When considering fractional exponents, such as $\alpha = 1.5$, the simulations exhibited heavier tails than for the classical diffusion order of $\alpha = 1.99$. However, increasing the domain size and refining the time step have proven to be effective in controlling the loss of mass, particularly for values of $\alpha$ close to 2. Furthermore, agent-based simulations provide a complementary perspective, illustrating the fact that fractional diffusion leads to a higher frequency of long-range particle jumps, an observation consistent with the statistical properties of stable distributions.
A major outcome of this study is the link between microscopic and macroscopic models, which indicates that fractional diffusion equations emerge as the macroscopic limit of a CTRW governed by heavy-tailed probability distributions. This connection bridges individual particle-based models with continuum descriptions, and further fortifies the role of fractional calculus in multiscale modelling.

Altogether, fractional calculus is increasingly being recognised as a powerful mathematical framework for capturing non-local effects, memory phenomena, and anomalous transport in complex systems. However, several research directions remain open for further exploration. 
An area of improvement can be found in the development of enhanced numerical methods: future efforts could focus on refining higher-order numerical schemes, and implementing adaptive mesh refinement strategies to increase accuracy. Another important direction would be given by the extension of these methods to higher-dimensional domains: many biological systems are better described in two- or three-dimensional (heterogeneous) environments. Extending the numerical framework to account for multidimensional dynamics would significantly broaden the scope and applicability of these models, particularly in the modelling of tissue structures, tumour progression, or cellular transport processes. Beyond dimensionality, incorporating reaction terms into the fractional diffusion equations would allow for a more comprehensive description of various biological (or other) processes, such as oncolytic virotherapy, with the goal of investigating how different treatment strategies affect tumour progression and how optimal therapeutic interventions can be designed. Such approaches would bring fractional models closer to real-world applications, where transport phenomena are coupled with reaction dynamics. 
Additionally---and even more importantly---validation against real biological data remains a necessary next step. This would be crucial in confirming the predictive capabilities of these models, confirming or even justifying the use of fractional calculus in the corresponding problems.

% %=========================== Thanks
% \section*{Acknowledgments} 
% We would like to acknowledge the assistance of volunteers in putting together this example manuscript and supplement. 

%=========================== Bibliography
\bibliographystyle{apalike}%{unsrtnat}%{apalike}%{abbrvnat}%{plainnat}%{siamplain}

%Abatangelo2023, Huang2014, Mao2017, Nochetto2015, Sheng2020,Cusimano2018,Anastasio1994,Li2015,
%Podlubny1999, Magin2006, Caffarelli2007, Henry2000, Metzler2000,  Baeumer2008

\bibliography{references}

\begin{thebibliography}{}

\bibitem[Abatangelo et~al., 2023]{Abatangelo2023}
Abatangelo, N., G\'omez-Castro, D., and V\'azquez, J.~L. (2023).
\newblock {Singular boundary behaviour and large solutions for fractional elliptic equations}.
\newblock {\em Journal of the London Mathematical Society}, 107(2):568--615.

\bibitem[Achleitner et~al., 2024]{Achleitner}
Achleitner, F., Akagi, G., Kuehn, C., Melenk, J.~M., Rademacher, J. D.~M., Soresina, C., and Yang, J. (2024).
\newblock {\em Fractional Dissipative PDEs}, pages 53--122.
\newblock Springer Nature Switzerland.

\bibitem[Anastasio, 1994]{Anastasio1994}
Anastasio, T.~J. (1994).
\newblock The fractional-order dynamics of brainstem vestibulo-oculomotor neurons.
\newblock {\em Biological Cybernetics}, 72(1):69--79.

\bibitem[Baeumer et~al., 2008]{Baeumer2008}
Baeumer, B., Kovács, M., and Meerschaert, M.~M. (2008).
\newblock Numerical solutions for fractional reaction–diffusion equations.
\newblock {\em Computers \& Mathematics with Applications}, 55(10):2212--2226.
\newblock Advanced Numerical Algorithms for Large-Scale Computations.

\bibitem[Bailo et~al., 2024]{Carrillo}
Bailo, R., Carrillo, J.~A., Fronzoni, S., and G\`omez-Castro, D. (2024).
\newblock {A finite-volume scheme for fractional diffusion on bounded domains}.
\newblock {\em European Journal of Applied Mathematics}, pages 1--21.

\bibitem[Benson et~al., 2000]{Benson2000}
Benson, D.~A., Wheatcraft, S.~W., and Meerschaert, M.~M. (2000).
\newblock Application of a fractional advection-dispersion equation.
\newblock {\em Water Resources Research}, 36(6):1403--1412.

\bibitem[Bezanson et~al., 2017]{Juliabezanson2017}
Bezanson, J., Edelman, A., Karpinski, S., and Shah, V.~B. (2017).
\newblock {Julia: A fresh approach to numerical computing}.
\newblock {\em SIAM review}, 59(1):65--98.

\bibitem[Caffarelli and Silvestre, 2007]{Caffarelli2007}
Caffarelli, L. and Silvestre, L. (2007).
\newblock An extension problem related to the fractional laplacian.
\newblock {\em Communications in Partial Differential Equations}, 32(7--9):1245--1260.

\bibitem[Cusimano et~al., 2018]{Cusimano2018}
Cusimano, N., del Teso, F., Gerardo-Giorda, L., and Pagnini, G. (2018).
\newblock {Discretizations of the spectral fractional Laplacian on general domains with Dirichlet, Neumann, and Robin boundary conditions}.
\newblock {\em SIAM Journal on Numerical Analysis}, 56(3):1243--1272.

\bibitem[Ding et~al., 2014]{scheme}
Ding, H., Li, C., and Chen, Y. (2014).
\newblock {High-Order Algorithms for Riesz Derivative and Their Applications}.
\newblock {\em Journal of Computational Physics}, 2014.

\bibitem[Gutleb and Carrillo, 2023]{Carrillo2023}
Gutleb, T.~S. and Carrillo, J.~A. (2023).
\newblock A static memory sparse spectral method for time-fractional pdes.
\newblock {\em Journal of Computational Physics}, 494:112522.

\bibitem[Henry and Wearne, 2000]{Henry2000}
Henry, B. and Wearne, S. (2000).
\newblock Fractional reaction–diffusion.
\newblock {\em Physica A: Statistical Mechanics and its Applications}, 276(3):448--455.

\bibitem[Huang and Oberman, 2014]{Huang2014}
Huang, Y. and Oberman, A. (2014).
\newblock {Numerical methods for the fractional Laplacian: A finite difference-quadrature approach}.
\newblock {\em SIAM Journal on Numerical Analysis}, 52(6):3056--3084.

\bibitem[Kelly et~al., 2019]{BCtwosided}
Kelly, J.~F., Sankaranarayanan, H., and Meerschaert, M.~M. (2019).
\newblock {Boundary conditions for two-sided fractional diffusion}.
\newblock {\em Journal of Computational Physics}, 376:1089--1107.

\bibitem[Kleinhans and Friedrich, 2007]{Kleinhans2007}
Kleinhans, D. and Friedrich, R. (2007).
\newblock Continuous-time random walks: Simulation of continuous trajectories.
\newblock {\em Physical Review E}, 76(6):061102.

\bibitem[Li and Zeng, 2015]{Li2015}
Li, C. and Zeng, F. (2015).
\newblock {\em Numerical Methods for Fractional Calculus}.
\newblock Chapman and Hall/CRC.

\bibitem[Magin, 2006]{Magin2006}
Magin, R.~L. (2006).
\newblock Fractional calculus in bioengineering.
\newblock {\em Critical Reviews in Biomedical Engineering}, 34(4):357--408.

\bibitem[Mainardi, 2010]{Mainardi2010}
Mainardi, F. (2010).
\newblock {\em Fractional Calculus and Waves in Linear Viscoelasticity: An Introduction to Mathematical Models}.
\newblock World Scientific.

\bibitem[Mainardi et~al., 2001]{analytic}
Mainardi, F., Luchko, Y., and Pagnini, G. (2001).
\newblock {The fundamental solution of the space-time fractional diffusion equation}.
\newblock {\em Fractional Calculus and Applied Analysis}, 4(2):153--192.

\bibitem[Mao and Shen, 2017]{Mao2017}
Mao, Z. and Shen, J. (2017).
\newblock {Hermite spectral methods for fractional PDEs in unbounded domains}.
\newblock {\em SIAM Journal on Scientific Computing}, 39(5):A1928--A1950.

\bibitem[M{\'e}ndez et~al., 2010]{reaction-transport}
M{\'e}ndez, V., Fedotov, S., and Horsthemke, W. (2010).
\newblock {\em {Reaction-transport systems: Mesoscopic foundations, fronts, and spatial instabilities}}.
\newblock Springer Series in Synergetics. Springer, Heidelberg.

\bibitem[Metzler and Klafter, 2000]{Metzler2000}
Metzler, R. and Klafter, J. (2000).
\newblock The random walk's guide to anomalous diffusion: a fractional dynamics approach.
\newblock {\em Physics Reports}, 339(1):1--77.

\bibitem[Nochetto et~al., 2015]{Nochetto2015}
Nochetto, R.~H., Ot\'arola, E., and Salgado, A.~J. (2015).
\newblock {A PDE approach to fractional diffusion in general domains: A priori error analysis}.
\newblock {\em Foundations of Computational Mathematics}, 15(3):733--791.

\bibitem[Nolan, 2020]{Nolan}
Nolan, J.~P. (2020).
\newblock {\em {Univariate Stable Distributions: Models for Heavy-Tailed Data}}.
\newblock Springer Series in Operations Research and Financial Engineering. Springer.

\bibitem[Omelchenko, 2010]{MultiStable}
Omelchenko, V. (2010).
\newblock {Elliptical Stable Distributions}.

\bibitem[Ortigueira, 2013]{FracTheory}
Ortigueira, M.~D. (2013).
\newblock {\em {Fractional Calculus for Scientists and Engineers}}.
\newblock Springer Publishing Company, Incorporated.

\bibitem[Podlubny, 1999]{Podlubny1999}
Podlubny, I. (1999).
\newblock {\em Fractional Differential Equations: An Introduction to Fractional Derivatives, Fractional Differential Equations, to Methods of Their Solution and Some of Their Applications}.
\newblock Academic Press.

\bibitem[Scalas et~al., 2004]{SubdiffFD}
Scalas, E., Gorenflo, R., and Mainardi, F. (2004).
\newblock {Uncoupled continuous-time random walks: Solution and limiting behavior of the master equation}.
\newblock {\em Phys. Rev. E}, 69:011107.

\bibitem[Sheng et~al., 2020]{Sheng2020}
Sheng, C., Shen, J., Tang, T., Wang, L.-L., and Yuan, H. (2020).
\newblock {Fast Fourier-like mapped Chebyshev spectral-Galerkin methods for PDEs with integral fractional Laplacian in unbounded domains}.
\newblock {\em SIAM Journal on Numerical Analysis}, 58(5):2435--2464.

\bibitem[\'Sl\k{e}zak, 2025]{github}
\'Sl\k{e}zak, J. (2025).
\newblock {StableDistributions.jl}.
\newblock \url{https://github.com/jaksle/StableDistributions.jl}.

\bibitem[Vieira et~al., 2023]{Motivation}
Vieira, L.~C., Costa, R.~S., and Valério, D. (2023).
\newblock {An Overview of Mathematical Modelling in Cancer Research: Fractional Calculus as Modelling Tool}.
\newblock {\em Fractal and Fractional}, 7(8).

\bibitem[West, 2016]{West2016}
West, B.~J. (2016).
\newblock {\em Fractional Calculus View of Complexity: Tomorrow's Science}.
\newblock CRC Press.

\bibitem[Xie and Fang, 2022]{XieFang}
Xie, C. and Fang, S. (2022).
\newblock {Efficient numerical methods for Riesz space-fractional diffusion equations with fractional Neumann boundary conditions}.
\newblock {\em Applied Numerical Mathematics}, 176:1--18.

\end{thebibliography}

\newpage

\end{document}